\documentclass[12pt,reqno]{article}
\usepackage{graphics}
\newcommand{\Om}{\Omega}
\newcommand{\Omr}{\Omega_r}
\newcommand{\Omz}{\Omega_0}

\renewcommand{\div}{{\rm div\, }}
\newcommand{\curl}{{\rm curl\, }}

\newcommand{\RR}{{\bf R}}

\newcommand{\EE}{{\cal E}_g}
\newcommand{\TT}{{\cal T}}
\newcommand{\SSS}{{\cal S}}

\newcommand{\Pf}{\noindent {\it Proof: \  }}
\newcommand{\QED}{\newline $\diamondsuit$}
\newcommand{\iint}{\int\!\!\!\!\int}

\newtheorem{thm}{Theorem}[section]
\newtheorem{lem}[thm]{Lemma}

\newtheorem{prop}[thm]{Proposition}
\newtheorem{rem}[thm]{Remark}

\newcommand{\be}{\begin{equation}}
\newcommand{\ee}{\end{equation}}
\newcommand{\bea}{\begin{eqnarray}}
\newcommand{\eea}{\end{eqnarray}}
\newcommand{\beann}{\begin{eqnarray*}}
\newcommand{\eeann}{\end{eqnarray*}}
\newcommand{\nnn}{\nonumber}

\begin{document}
\title{
Minimizers of the Lawrence--Doniach energy\\
in the small-coupling limit: \\
finite width samples in a parallel field}
\author{{\Large S. Alama\footnote{Dept. of Mathematics and Statistics,
McMaster Univ., Hamilton, Ontario, Canada L8S 4K1.  Supported
by an NSERC Research Grant.}, \ A.J. Berlinsky%
\footnote{Dept. of Physics and Astronomy, McMaster Univ., Hamilton,
Ontario, Canada L8S 4K1.  Supported by an NSERC Research Grant.}, \
\&  L. Bronsard${}^*$}}
\thispagestyle{empty}
\maketitle

\begin{abstract}
In this paper we study the Lawrence--Doniach model for layered
superconductors, for a sample with finite width subjected to a
magnetic field parallel to the superconducting layers.  We
provide a rigorous analysis of the energy minimizers in the limit
as the coupling between adjacent superconducting layers
tends to zero. We identify a unique global minimizer of the Gibbs free
energy in this regime (``vortex planes''),
and reveal a sequence of first-order phase
transitions by which Josephson vortices are nucleated via
the boundary.  The small coupling limit is studied via degenerate
perturbation theory based on a Lyapunov--Schmidt decomposition
which reduces the Lawrence-Doniach system to a finite-dimensional
variational problem. Finally, a lower bound on the radius of validity of
the perturbation expansion (in terms of various parameters
appearing in the model) is obtained.

\end{abstract}

\newpage

\baselineskip=17pt

\section{Introduction}

In 1971 
Lawrence and Doniach \cite{LD} introduced a 
Ginzburg--Landau type model for superconducting materials
with a planar layered structure.   In this model, the
superconductor occupies an array of parallel sheets
with insulating material acting as a buffer between the
sheets.  While this model was originally proposed to
study layered structures artificially produced by
successively deposing thin planar sheets of superconducting metal
with organic separators, it has received renewed attention
due to the discovery of  high
temperature superconductors.  Indeed,  nearly all of these high-$T_c$
materials are crystals with
 a distinctly layered structure, consisting of
copper oxide superconducting planes stacked with
intervening insulating (or weakly superconducting) planes.

\smallskip

In this paper we will consider the case of a layered superconductor
in a uniform magnetic field imposed  parallel to the
superconducting planes.  We assume that there
are a finite number of superconducting sheets, each parallel to the
$xy$-plane, with uniform separation $p$.
We assume that the external magnetic field is applied along the $y$-direction,
$\vec H=H\, \hat y$.  
We will take the planes to be of infinite
extent in the $y$-direction, and assume the
local magnetic field will be everywhere independent of
$y$ and point in the $y$-direction, 
$$    \vec h(x,y,z)=h(x,z)\, \hat y.   $$
The vector potential $\vec A$ may then be chosen to lie
in the $xz$-plane,
$$  \vec A(x,y,z)= A_x(x,z)\, \hat x \, + A_z(x,z)\, \hat z,
\qquad  \vec h = \curl \vec A = 
    \left(  {\partial A_x\over \partial z} -
         {\partial A_z\over \partial x} \right) \, \hat y.  $$
We assume that the sample has fixed width $2L$ in the
$x$-direction, and hence the superconducting sheets are
described by the stack of parallel planar strips
$$  \Sigma_N:\  -L\le x\le L, \quad  -\infty<y<\infty, \quad
   z=z_n:=np, \ n=0,1,\dots,N.  $$

Since each sheet $\Sigma_N$ is superconducting
it carries  a (complex-valued) order parameter
$\psi_n(x)$, $n=0,\dots,N$.  We choose units in such a way that
$|\psi_n|=1$ represents a purely superconducting
state.  The Lawrence-Doniach model is then formulated in terms
 of the following Gibbs free energy functional:
\beann
{\cal G}_r (\psi_n,\vec A)&=& {H_c^2\over 4\pi} \left\{
   p\sum_{n=0}^N
    \int_{-L}^{L}  \left[
        {1\over \kappa^2}
           \left|\left( {d\over dx} -iA_x\right)\psi_n  \right|^2 +
       \frac12 (|\psi_n|^2-1)^2  \right]\, dx 
         \right. \\   
             &&\qquad + {r\over 2}\, p\sum_{n=1}^N
                  \int_{-L}^L  \left| \psi_n - 
                     \psi_{n-1}\exp\left(i\int_{z_{n-1}}^{z_n} A_z(x,s)\, ds\right)
                          \right|^2
                  \, dx \\  
            &&  \left.
            \qquad   +\ {1\over  \kappa^2}
                \int_{-L}^{L}\int_0^{Np} \left(  {\partial A_x\over \partial z} 
                    - {\partial A_z\over \partial x}
      - H\right)^2 \, dz\, dx \right\},
  \eeann
where $r$ is the
{\it interlayer coupling parameter}
(or {\it Josephson coupling parameter}.)
We have chosen units such that the
in-plane penetration depth $\lambda_{ab}=1$,
$\kappa=\lambda_{ab}/\xi_{ab}$ is the Ginzburg--Landau parameter, 
($\lambda_{ab},\xi_{ab}$ are the in-plane penetration depth
and correlation length, respectively,) 
and the magnetic fields
are measured in units of $H_c/\kappa$,
where  $H_c$ is the thermodynamic critical field. (See \cite{Tinkham}.)

The coupling between the superconducting planes given
by the second sum in ${\cal G}_r$ simulates the Josephson
effect, by which superconducting electrons travel
from one superconducting sheet to another by quantum mechanical
tunnelling.  We will see this explicitly in the Euler--Lagrange
equations, where the currents in the gaps between planes will
be determined by the sine of the gauge-invariant phase difference.
The interlayer coupling parameter $r$ gives the strength of
the Josephson effect.  In our units,
$$  r={2\over \lambda_J^2 \kappa^2 p^2}  $$
where $\lambda_J$ is the 
(non-dimensional) {\it Josephson penetration depth.}
In an anisotropic Ginzburg--Landau model, $\lambda_J^2$
gives the effective mass ratio which determines the degree
of anisotropy.  For highly anisotropic superconductors
$\lambda_J$ is very large, and hence when $\kappa p\sim 1$
we may treat $r$ as a small parameter.

\smallskip

Due to the layered structure one expects these materials
to be highly anisotropic.  A first attempt to model layered
superconductors is by an {\it anisotropic Ginzburg--Landau}
model, which treats the sample as a three-dimensional
solid with anisotropic material parameters.  
For certain materials  and temperatures close
to the critical temperature $T_c$ this approximation
seems valid, but for the most anisotropic superconductors
the anisotropic Ginzburg--Landau model does not give a
good qualitative or quantitative description of experimental
observations.  For example, when the sample is subjected to a sufficiently
strong magnetic field oriented {\it parallel}\, to the superconducting
planes  Kes,  Aarts,  Vinokur, and van der Beek \cite{Kes} 
observe a transition betwen ``three-dimensional'' behavior 
(governed by the anisotropic Ginzburg--Landau model)
and ``two-dimensional'' behavior at a critical temperature $T_{C0}$
below $T_C$.  In the two-dimensional regime the superconducting
planes decouple and the applied magnetic
field penetrates completely
between the planes, virtually unscreened by the superconductor.
Despite the penetration of the field, superconductivity within the 
planes is not destroyed even in very strong applied fields.
This ``magnetically transparent'' state is inconsistent with
the anisotropic Ginzburg--Landau model, where  magnetic
fields of moderate intensity are largely expelled from the bulk 
except for an array of isolated
vortices (the ``Abrikosov lattice''.)
In addition the Ginzburg--Landau model predicts the
breakdown of superconductivity when the applied field penetrates
the material completely, 
which occurs when the field exceeds a critical value $H_{c2}$.

\smallskip

We note that Chapman, Du, \& Gunzburger \cite{CDG} have proven
that solutions of the Lawrence--Doniach model converge to
solutions of the anisotropic Ginzburg--Landau model (and
in particular the convergence of energy minimizers) under the
limit $p\to 0$ with $\kappa,\lambda_J$ fixed.  This limit does not 
correspond to our ``two-dimensional'' regime, since it
would send $r\to\infty$, corresponding to a strong
coupling between adjacent superconducting layers.
Indeed, we observe that
the non-dimensional separation distance $p$ of our model is
related to the (dimensionally dependent) physical separation
$\bar p$ via $p=\bar p/\lambda_{ab}$.  Since 
$\lambda_{ab}(T)\to\infty$ as $T\to T_c$, the limit
$T\to T_c$ effectively corresponds to $p\to 0$ 
(with $\kappa,\lambda_J,\bar p$ fixed) in our units,
and therefore the Chapman, Du, \& Gunzburger
limit can be interpreted as letting $T\to T_c$.
This is consistent with the observed  ``dimensional
crossover'' to the anisotropic Ginzburg--Landau regime
for temperatures near $T_c$.

\medskip

In this paper we will study the minimizers (and low-energy
solutions) of the Lawrence--Doniach system for $r$ near zero
and thereby analyse the structure of the resulting ``transparent state.''
The crucial observation is that when $r=0$ the planes decouple,
and the energy may be minimized explicitly by solving
simple {\it first order} equations.    
Even after gauge symmetries
have been removed the  $r=0$ problem exhibits an additional
symmetry, corresponding to an $N$-dimensional torus action (where
$N+1$ denotes the number of superconducting planes,)  
and thus the minimization problem at $r=0$ degenerates on
 a finite dimensional manifold in function space.
We can think of the $r=0$ problem in analogy with the
self-dual point of the Ginzburg--Landau model ($\kappa=1/\sqrt{2}$),
where minimizers satisfy a first-order Bogomolnyi system
(in addition to the usual second order Ginzburg--Landau equations,)
and the same minimum energy is obtained by any
configuration of vortices.

When $r\neq 0$ this symmetry is broken and 
 a Lyapunov--Schmidt decomposition  reduces
the problem of finding solutions with $r\simeq 0$ to a
finite dimensional variational problem on the degenerate manifold.
The minimum value of energy is
$O(r)$, and we indeed recover the ``transparent state''
observed in experiments.  The local magnetic field
$h(x,z)=H + O(r)$ inside the sample, and superconductivity
is hardly affected in each plane, $|\psi_n(x)|=1-O(r)$.
In particular the order parameters are never zero:  ``vortices''
correspond to local maxima of the local magnetic field, and
lie between the layers.  In the physics literature these
are referred to as {\it Josephson vortices,} as opposed to
the
Abrikosov vortices typically observed in the Ginzburg--Landau
theory.

The finite dimensional reduced problem
may be solved explicitly to determine the exact
geometry of the Josephson vortex lattice.
For a finite sample in $x,z$ the minima of energy form
``vortex plane''
configurations, in which the local magnetic field is uniform
in $z$.  The vortices are not separated, but line up vertically
at the local maxima of $h(x,z)=h(x)$.  (See figure 1.)
As the external field
$H$ is increased vortices are nucleated at the edges, by
a first-order phase transition.  Formal asymptotic
expansions for these solutions were calculated 
by Theorodakis \cite{Th}, and Kuplevakhsky \cite{K}
claimed that they were the only solutions of the 
Lawrence--Doniach system.  On the contrary, we find that there are
exactly $2^N$ solutions with energy $O(r)$:  two represent
vortex planes (one stable and the other
unstable), and the others (unstable) lattices of
various geometries.  (See Theorem~\ref{Thm1}.)

Again, we note the distinction with the Ginzburg--Landau model:
the geometry of the Abrikosov lattice was determined by
numerical comparison of a finite number of possible lattice
geometries.  For the Lawrence--Doniach model in the small
coupling limit we are able to identify the absolute minimizer
(and all low-energy solutions of the Euler--Lagrange equations)
explicitly and rigorously.  This is a direct benefit of the
discrete nature of the model.

 The basic idea that an infinite dimensional
variational problem is actually governed by a finite dimensional
one in some parameter limit is a common one in analysis.  Indeed,
it appears in such diverse contexts as the location of 
Ginzburg--Landau vortices as $\kappa\to\infty$ (Bethuel, Brezis,
\& H\'elein \cite{BBH}, Bethuel \& Rivi\'ere \cite{BR},)
spike-layer solutions  (for example, 
Li \& Nirenberg \cite{LN}, Gui \cite{Gui}, or Wei \cite{W}), and blow-up
for critical exponent problems (for example in Bahri, Li, \& Rey
\cite{BLR}, or Rey \cite{Rey}.)   These examples are of {\it singular}
perturbation problems, though.  The transparent state arises as a {\it
degenerate regular} perturbation of the $r=0$ problem, and hence it is more
closely related to the work of Ambrosetti, Coti-Zelati, \& Ekeland
\cite{ACE} and Ambrosetti \& Badiale \cite{AB} on homoclinic
solutions of Hamiltonian systems and the Poincar\'e-Melnikov
functions.

\medskip

Analytically our results are unambiguous:  for any choice of the
other parameters ($L$, $\kappa$, $N$, $H$, $p$)  we can
choose $r$ sufficiently small so that the vortex planes configuration
minimizes the free energy.  However, in a real superconductor
$r$ is not infinitessimally small, which raises the question
of how the interval of validity of the $r$-expansion is affected
by the values of the other parameters in the problem.
For example, in the order $r$ term in the expansion
of the solution (see (\ref{Phi}), for example)
we observe `` secular'' terms, that is factors which become
large without bound as the length of the interval increases,
and in general remark that the coefficients increase with
$L,\kappa$ and decrease with $H$.
We address this question in section~5, where we produce a
lower bound for the interval of validity as a function of
$N$, $L$, $\kappa$, and $H$.  We discover that this interval
is independent of $N$ and indeed increases with decreasing $L$, $\kappa$
and increasing $H$.  For the high-$T_c$ superconductors,
$\kappa$ is large and typical macroscopic sample widths $L$
are very large compared with the in-plane penetration depth
$\lambda_{ab}$.
This indicates that our analysis may be more applicable to
experiments with highly anisotropic organic or synthetic
multi-layer superconductors, where the material parameters are
significantly different from the high-$T_c$ crystals.

\smallskip

In a subsequent paper \cite{ABB2} we address the question of
minimizers of the Lawrence--Doniach energy in very large samples,
by considering {\it periodic} solutions in an
 infinitely wide sample.  In fact,
in the periodic case we find a different geometry
for the energy minimizing configuration!  The periodic solution with the
least energy is  a lattice with
period two in $n$, forming a diamond pattern of Josephson
vortices, proposed by
Bulaevski\u \i  \ \& Clem \cite{BC}.
In the periodic case, the role of $L$ is now taken by the
period, which for minimizers will decrease as the
applied field $H$ increases.  Therefore the small $r$
expansion will have a large range of validity in
sufficiently large fields $H$, and the result should
better describe experiments with high-$T_c$ superconductors.

\bigskip

\noindent
{\bf Acknowledgement.}\
The authors thank A. Bahri for his helpful suggestions
for determining the region of validity of the perturbation
methods.
SA and LB wish to thank R. V. Kohn and P. A. Deift for inviting
them to spend a semester at the Courant Institute in Fall 1999,
and for the interest and encouragement they have shown 
throughout the years.

\setcounter{thm}{0}

\subsection{Variational setting}
We begin with the following basic energy estimate, which
legitimizes the simplification $|\psi_n(x)|\simeq 1$
for small $r$.  Indeed, in the physics literature the approximation
$|\psi_n|\sim 1$ is widely assumed to hold:
see for example Bulaevskii \cite{Bul}, Clem \& Coffey \cite{CC}, or
Bulaevskii \& Clem \cite{BC}.
In particular, we see immediately that in the regime $r\simeq 0$
there will be no ``Abrikosov vortices'' in the sense that the
order parameter is never zero.
\begin{prop}\label{basicenergy}
Let $D:=[-L,L]\times[0,Np]$.
For any $r\ge 0$ we have
$$  \inf \{ {4\pi\over H_c^2}{\cal G}_r (\psi_n,\vec A): \ 
\psi_n\in H^1([-L,L]),
  \  \vec A\in H^1(D;{\bf R}^2)\} \, \le \, 2Np(L+{1\over pH})r.  $$
Moreover, there exist constants $r_0>0$ and $C=C(N,L,\kappa,H,p)>0$ such that
for all $r\in [0,r_0)$ 
 the minimum is attained by $(\psi_n,\vec A)$
with $|\psi_n(x)|\ge 1-Cr^{1/2}>0$ for all $x\in [-L,L]$.
\end{prop}

\noindent
\Pf
We choose a test configuration, $\psi_n(x)=\exp\{inpHx\}$, $A=(Hz,0)$
to estimate the free energy,
$$  {4\pi\over H_c^2}{\cal G}_r(\psi_n,A) = rp \sum_{n=1}^N  \int_{-L}^L
    (1-cos(pHx))\, dx \le 2N\left(L+{1\over pH}\right)pr.  $$
The fact that the minimum is attained can be easily proven once
an appropriate choice of gauge has been made: see \cite{CDG} for
details.

It remains to show that the order parameters lie near the unit
circle.
We recall the ``Diamagnetic inequality'' (see for example
p. 174 of \cite{LiebLoss}), which states 
$$  \left|  \nabla |f| \right|(x) \le |(\nabla - iA)f|(x)
\quad\mbox{(almost everywhere)}  $$
for every $f\in L^2_{loc}$ with $(\nabla - iA)f\in L^2_{loc}$.
Using the elementary inequality $(1-|\psi_n|)^2\le (1-|\psi_n|^2)^2$,
and the energy bound we obtain 
\beann
  \sum_{n=0}^N \| 1-|\psi_n| \|_{H^1}^2
    & \le & \sum_{n=0}^N \int_{-L}^L  \left( (1-|\psi_n|^2)^2
        + {1\over \kappa^2} \left| {d\over dx} |\psi_n| \right|^2 \right)\, dx\\
         &\le & \sum_{n=0}^N \int_{-L}^L  \left( (1-|\psi_n|^2)^2
        + {1\over \kappa^2} |({d\over dx} -iA_x)\psi_n|^2 \right)\, dx \\
        & \le &  {4\pi\over H_c^2} {\cal G}_r(\psi_n,A) 
            \le 2N\left(L+{1\over Hp}\right)pr.
 \eeann
Here we choose an equivalent norm for $H^1([-L,L])$,
$$  \| f\|_{H^1} :=\sqrt{\int_{-L}^L 
     \left( |f|^2 + {1\over\kappa^2} |f'|^2 \right)\, dx }.
$$
By the Sobolev embedding we have for each $n=0,\dots,N$
that $|\psi_n(x)|\ge 1- Cr^{1/2}$ with constant $C$ depending
on $N, L, p, \kappa,H$.
\QED

The above proposition suggests the use of polar coordinates for $\psi_n$
in order to deal more directly with the phase of the order parameter,
which plays the essential role in Josephson coupling.
We define $f_n,\phi_n$ via $\psi_n(x)=f_n(x)\exp (i\phi_n(x))$,
and note that $\phi_n$ is well-defined only up to an additive
integer multiple of $2\pi$.   We then define our new free energy
in terms of the variables $(f_n,\phi_n,\vec A)$ to coincide with
$ {4\pi\over H_c^2}{\cal G}_r$, that is:
\beann
\Omega_r (f_n,\phi_n,\vec A)&=& p\sum_{n=0}^N 
    \int_{-L}^{L}\left[
       \frac 12  (f_n^2-1)^2 + {1\over\kappa^2} (f'_n)^2
            + {1\over\kappa^2}(\phi'_n- A_x(x,z_n))^2 f_n^2
             \right]\, dx \\   
             &&\qquad + {r\over 2} p\sum_{n=1}^N   \int_{-L}^L
                 \left( f_n^2 +f_{n-1}^2 
              -2 f_n f_{n-1}\cos(\Phi_{n,n-1} )\right) \, dx \\  
            &&\qquad   +\ {1\over  \kappa^2}
                \int_{-L}^{L}\int_0^{Np} \left(  {\partial A_x\over \partial z} 
                    - {\partial A_z\over \partial x}
      - H\right)^2 \, dz\, dx ,
  \eeann
where 
$$\Phi_{n,n-1}(x):= \phi_n-\phi_{n-1}-
                \int_{z_{n-1}}^{z_n} A_z(x,z)\, dz, $$
is the {\it gauge-invariant phase difference}, 
$z_n=np$ and
$r={2\over \lambda_J^2 \kappa^2 p^2}$.
When $f_n=|\psi_n|$ is bounded away from zero
the condition $\psi_n\in H^1([-L,L])$ is equivalent to 
both $f_n,\phi_n\in H^1([-L,L])$.

We first define a base space, in which $\Omega_r$ will be
a smooth functional:
$$
{\cal E}:=
\left\{  \begin{array}{c}
(f_n,\phi_n,\vec A): \ f_n\in H^1([-L,L]), \phi_n\in H^1([-L,L]), \ n=0,\dots, N \\
   \vec A=(A_x,A_z)\in H^1(D,{\bf R}^2)
\end{array}  \right\}.
$$
We remark that we should really work in the convex subset of
${\cal E}$ with $f_n\ge 0$ for all $n=0,\dots,N$, but  
Lemma~\ref{basicenergy} already guarantees that the solutions
we will find will have $f_n\sim 1$.  Furthermore, Proposition~4.5
of \cite{CDG} asserts that $f_n(x)=|\psi_n(x)|\le 1$ for any solution
of the Lawrence--Doniach equations.

We note that $\Omega_r$ is a smooth ($C^\infty$) functional
on $\cal{E}$, and that variation with respect to each of its arguments
gives the Lawrence--Doniach system.  First, we denote by
\be
\label{supervelocity}
V_n:=(\phi'_n- A_x(x,z_n)),
\ee
the supercurrent velocity.
Then,  variation of $\Omr$
with respect to $f_n$ (for each $n=0,\dots,N$) yields
\be\label{feqn}
 -{1\over \kappa^2} f''_n + (f_n^2-1)f_n 
    + {1\over\kappa^2}V_n^2 f_n =
    \cases{
          {r\over 2}\left(
       f_{n-1}\cos\Phi_{n,n-1}
      +f_{n+1}\cos\Phi_{n+1,n}-2f_n\right), & $n \neq 0,N$; \cr
  &\cr
     {r\over 2}\left(f_{1}\cos\Phi_{1,0}-f_0\right), & $n=0$; \cr
  &\cr
      {r\over 2}\left(
       f_{N-1}\cos\Phi_{N,N-1}
         -f_N\right),  & $n=N$,\cr}
\ee
with boundary condition $f'_n(\pm L)=0$, $n=0,\dots, N$.

Variation with respect to $A_x$ produces the following
equation in weak form:
\bea
\label{hz}
&&{\partial h\over \partial z}(x,z)=0 \qquad \mbox{for all 
   $(x,z)\in (-L,L)\times (z_{n-1}z_n)$,
$n=0,\dots,N$,}
\\  \label{jump}
&& h(x,z_n+)-h(x,z_n-) = - p f_n^2(x) (\phi'_n - A_x(x,z_n)),
\quad n=1,\dots,N-1, \\
\nnn
&& h(x,z_N-)=p f_N^2(x) (\phi'_N - A_x(x,z_N), \quad
   h(x,0+)=-p f_0^2(x) (\phi'_0 - A_x(x,0).
 \eea
 In other words, $h$ is independent of $z$ away from the superconducting
planes, and supercurrents in the SC planes
create jump discontinuities.
The dependence of $h(x,z)$ on $x$ in each gap is determined by (compactly
supported) variations of $A_z$: 
\bea\label{hx}
{\partial h\over \partial x} =
     {r\kappa^2 p \over 2} f_n(x) f_{n-1}(x) \sin\Phi_{n,n-1}(x),
           &&\mbox{ if $-L\le x\le L$ and}\\
           \nnn
       &&  z_{n-1}<z<z_n, n=1,\dots,N;
\eea
with boundary condition $h(\pm L,z)=H$.
It is therefore natural to define
\bea\label{heqn1}
&&  h(x,z)=h^{(n)}(x), \quad \mbox{when $-L\le x\le L$, $z_{n-1}<z<z_n$,
   $ n=1,\dots,N$,}
\eea
with $h^{(n)}$ determined by the 
ordinary differential equations (\ref{hx}) together with the boundary
condition $h^{(n)}(\pm L)=H$.

Variation with respect to $\phi_n$ produces the following
current--conservation laws:  
\be\label{current}
{1\over \kappa^2}{d\over dx}\left( f_n^2(\phi'_n - A_x(x,z_n))\right) =
   \cases{
       {r\over 2}[f_n  f_{n-1}\sin\Phi_{n,n-1}
       -f_{n+1}f_n \sin\Phi_{n+1,n}], & $n=1,\dots,N-1$;\cr &\cr
    -{r\over 2} f_1 f_0 \sin\Phi_{1,0}, \
         &  $n=0$; \cr &\cr
     {r\over 2} f_{N} f_{N-1} \sin\Phi_{N,N-1},
          & $n=N$,   }
 \ee
with boundary condition 
$$  f_n^2(x)(\phi'_n(x) - A_x(x,z_n)) = 0, \quad x=\pm L, \  n=0,\dots,N,$$
which expresses the physical fact that current should not
flow past the edge of the material.
 Note that (\ref{current})
 does not yield any new information, since it can be obtained
by differentiating the jump condition (\ref{jump}) and
substituting from (\ref{hx}).  This is not surprising, since
gauge invariance implies a nontrivial relationship between
the $\phi_n$ and  $\vec A$.  Indeed,
denoting  the supercurrents in the planes by
\be
\label{supercurrent}
j_x^{(n)}:=V_n f_n^2,
\ee
and the Josephson current in between the $n-1$ and the $n$th planes by
\be
\label{josephsoncurrent}
j_z^{(n)}:= {r\over 2}\kappa^2p f_nf_{n-1} \sin
\Phi_{n, n-1},
\ee
equation (\ref{current}) is a semi-discrete version of
the classical continuity equation $\div \vec j =0$, and 
gives the conservation law  corresponding to the $U(1)$
gauge invariance in accordance with Noether's
Theorem.

\smallskip  
   
   We are most interested in the gauge invariant ``observable''
quantities, which enter directly into the free energy:
the density of superconducting eletrons
$f_n$, the supercurrent velocity
$V_n$, the gauge-invariant phase difference
$\Phi_{n,n-1}$, and the local magnetic field $h(x,z)$.
A very useful formula for $\Phi_{n,n-1}(x)$ can be
obtained by applying Stokes' Theorem in the rectangle
$R=(0,x)\times (z_{n-1},z_n)$:
\be\label{Phieqn}
\Phi_{n,n-1}(x) = \int_0^x (V_n-V_{n-1})\, d\bar x
   + p \int_0^x h^{(n)}(\bar x) \, d\bar x + \Phi_{n,n-1}(0)
   \quad n=1,\dots, N.
   \ee

From these equations we easily verify the smoothness of observables
associated to 
weak solutions of the Lawrence--Doniach system:
\begin{prop}\label{regularity}
Suppose $(f_n,\phi_n,\vec A)\in {\cal E}$ are critical points
of $\Omr$.  Then $f_n, V_n, \Phi_{n,n-1}\in C^\infty([-L,L])$, 
and $h\in C^\infty([-L,L]\times(z_{n-1},z_n))$, $n=1,\dots,N$.
\end{prop}
Of course, the regularity of the non-gauge-invariant quantities
$\phi_n$ and $\vec A$ depends on the choice of gauge.
\begin{rem} \rm
Some authors have (correctly) pointed out that 
it is not physically correct to impose the external field via
a Dirichlet condition $h=H$ on the boundary $\partial D$
of the sample.  A more appropriate model for the effect of an external
field is obtained by placing the superconductor $D$ in a larger 
region $\tilde D\supset D$, (with $\tilde D=\RR^2$
possibly) and including the field energy in $\tilde D\setminus D$
in the calculations of the free energy,
$$  \tilde\Omr(f_n,\phi_n,\vec A):= \Omr(f_n,\phi_n,\vec A) +
\iint_{\tilde D\setminus D}  (\curl \vec A -H )^2\, dx\, dz.
$$
For example, our sample $D$ may be lying at the center of
a long cylindrical solenoid of large radius (whose interior
cross-section is a large disk $\tilde D$.)
All of the preceding analysis can then be carried through for
$\tilde\Omr$ by choosing an appropriate space (see Rubinstein
\& Schatzman \cite{RS} for a discussion
of the correct setting when $\tilde D=\RR^2$,)
but it is easy to verify that the Euler-Lagrange equations
yield  $h(x,z)\equiv H$ for $(x,z)\in \tilde D\setminus D$.
This is due to the two-dimensional ansatz:  
$\nabla\times( g(x,z)\hat y) = 0$
in a domain
implies that $g$ is constant there.  Consequently we may
use the simpler form of the energy $\Omr$ with no loss of
generality or of physical relevance.
\end{rem}
 
There is a large degree of degeneracy of $\Omega_r$
in $\cal{E}$ due to the gauge invariance:  if $\chi\in H^2([-L,L]\times
[0,Np]))$, and
$$ \hat f_n=f_n,\qquad  \hat \phi_n(x) = \phi_n(x) - \chi(x,z_n), \qquad
  \hat{A}= \vec A - \nabla \chi,  $$
then $\Omr(\hat f_n, \hat \phi_n, \hat{A})=\Omr(f_n,\phi_n,\vec A)$. 
As usual, we eliminate this troublesome degeneracy by fixing a gauge.
The most convenient choice is the {\it Coulomb gauge}, which allows
us to control the $H^1$ norm of the vector potential $\vec A$ by
its curl.  We define a subspace of $\cal{E}$ to incorporate this choice
of gauge,
$$
\EE:= \{ (f_n,\phi_n,\vec A)\in {\cal E}: \ \int_{-L}^L\phi_0(x)\, dx=0, \
   \div{\vec A}=0  \  \mbox{in $D$, and} \
        \vec A\cdot \vec n = 0\ \mbox{on $\partial D$.}\}
$$
This choice is made with no loss of generality:
\begin{lem}\label{coulombfinie}
For every $(f_n,\phi_n,\vec A)\in {\cal E}$, there exists
$\chi\in H^2(D)$ so that
$(f_n,\phi_n-\chi(\cdot,z_n),\vec A-\nabla\chi)\in \EE.$
\end{lem}
We also have:
\begin{lem}\label{control}
For every $(f_n,\phi_n,\vec A)\in \EE$
\be\label{H1}
  \|\vec A\|_{H^1(D)}^2 \le  C_0 \|\curl \vec A\|_{L^2(D)}^2, 
\ee
where
\be\label{C_0}
     C_0  =   2\left[ 
         1 + {4\over \pi^2} {L^2 N^2 p^2\over N^2 p^2 + 4L^2}\right]^2.
\ee
\end{lem}

\noindent
{\it Proof of Lemma~\ref{coulombfinie}:}\
Given $(f_n,\phi_n,\vec A)\in {\cal E}$ we may solve
the linear Neumann problem,
$$ \left\{  \begin{array}{l} 
\Delta \chi = \div \vec A \quad \mbox{for $x\in D$}, \\
{\partial\chi\over\partial \vec n} = \vec A\cdot\vec n
 \quad \mbox{for $x\in\partial D$.}
\end{array}
\right.
$$
By the divergence theorem there exists a unique
(up to constants)  weak solution, $\chi\in H^1(D)$.
Since $\vec A\in H^1$ we have $\chi\in H^2_{loc}(D)$ by
standard regularity theory.  Since the domain $D$
is polygonal we must be more careful to determine
the regularity at the corners, but by Theorem~1.5.2.4
of Grisvard \cite{Grisvard} we may find an
$H^2(D)$ function whose normal derivative
coincides with  $\vec A\cdot\vec n$ on $\partial D$.
Then
Theorem~4.3.1.4 of \cite{Grisvard} provides global
regularity, $\chi\in H^2(D)$.  By subtracting a constant
from $\chi$ such that $\phi_0(x)-\chi(x,0)$ has average
zero we achieve the desired gauge change. 
\QED

\noindent
{\it Proof of Lemma~\ref{control}:}\
Assume $(f_n,\phi_n,\vec A)\in \EE$.  We now solve
the Dirichlet problem,
\be
\left\{  \begin{array}{l} 
\Delta \eta = \curl \vec A \quad \mbox{for $x\in D$}, \\
\nnn
\eta=0
 \quad \mbox{for $x\in\partial D$.}
\end{array}
\right.
\ee
By Theorem~4.3.1.4 of \cite{Grisvard}
the unique solution $\eta\in H^2(D)$, and 
$\tilde A=(\partial_z\eta, -\partial_x\eta)\in H^1(D;\RR^2)$
with $\div(\vec A- \tilde A)=0$, $\curl (\vec A-\tilde A)=0$, and
$(\vec A-\tilde A)\cdot\vec n=0$ on $\partial D$.
Hence $\vec A=\tilde A$, and the explicit solution 
of the Dirichlet problem in a rectangle (via Fourier analysis)
provides the constant $C_0>0$ such that:
$$  \|\vec A\|^2_{H^1(D)}
     \le \|\eta\|^2_{H^2(D)} 
         \le C_0 \|\Delta \eta\|^2_{L^2(D)}
      =   C_0\|\curl \vec A\|^2_{L^2(D)}.  $$
\QED

\section{Minimization at $r=0$.}
\setcounter{thm}{0}

When $r=0$, the superconducting planes decouple, and
we may solve the minimization problem explicitly.  The solution of the
problem $r\sim 0$ will require some detailed second-order
information on the minimizers at $r=0$, so we
first establish a functional analytic setting for the
equations.  We exploit the Hilbert
manifold structure of $\EE$, and regard the first variation of
energy as elements of the tangent Hilbert space
\beann      E=  T\EE &=&\left\{ (u_n,v_n,\vec a): \
   u_n,v_n\in H^1([-L,L]),\ \int_{-L}^L v_0\, dx=0,\right.  \\
  &&\quad \left.
    \vec a\in H^1(D), \
  \div\vec a=0\ \mbox{in $D$},\  \vec a\cdot\vec n=0 \
\mbox{on $\partial D$.} \right\}, 
\eeann
 and the second variation as a 
 self-adjoint operator on $E$.
We first introduce an equivalent inner product:  for 
$(u_n,v_n,\vec a),(U_n,V_n,\vec A)\in E$, let
\beann   
  &&  \qquad\qquad\left\langle 
       (u_n,v_n,\vec a), (U_n,V_n,\vec A)\right\rangle :=
\\
 & &  p\sum_{n=0}^N \int_{-L}^L \left\{ 2u_n U_n +
      {1\over\kappa^2}  [u'_n U'_n + v'_n V'_n] + v_n V_n\right\}\, dx 
        +  {1\over\kappa^2}\iint_D \curl \vec a \cdot \curl\vec A \, dx\, dz.
   \eeann
Then we define the gradient of $\Omr(f_n,\phi_n,\vec A)$ via
$$  \left\langle
   \nabla\Omr(f_n,\phi_n,\vec A) , (u_n,v_n,\vec a) 
       \right\rangle = D\Omr(f_n,\phi_n,\vec A)[u_n,v_n,\vec a].  $$

\begin{prop}\label{r=0}
With $r=0$,
\begin{enumerate}
\item[(a)] $\inf \{\Om_{0}(f_n,\phi_n,\vec A): \
   (f_n,\phi_n,\vec A)\in \EE\} \, =\, 0.$  The minimum value
   is attained, and the set of all minimizers coincides with
the hyperplane $\SSS$, defined by the set of
$(f_n,\phi_n,\vec A)\in \EE$ such that: \quad $f_n\equiv 1$;
\be\label{phi0}
     \phi_n(x)=\alpha_n +\int_0^x A_x(s,z_n)\, ds, 
        \quad \alpha_0=0, \quad \alpha_n\in {\bf R}, \ n=1,\dots,N; 
\ee
\be \label{A0}
             \vec A=(\eta_z, -\eta_x), \quad \mbox{where
                 $\Delta \eta =H$ in $D$, 
                    and  $\eta|_{\partial D}=0.$}
\ee
     In particular,  the
gauge invariant phase difference is given by $\Phi^0_{n,n-1}=\delta_n
     +Hpx$, where $\delta_n:= \alpha_n-\alpha_{n-1}$, $n=1,\dots,N$.
\item[(b)]  For any element $s=(f_n^0,\phi_n^0,\vec A^0)\in \SSS$, the
     linearized operator
    $D^2\Om_0(s): \ E\to E$ defines a Fredholm operator with index zero.
      Moreover,
 \be\label{nullspace}
    T_s\SSS=\ker D^2\Om_0(s)\simeq \RR^N. \ee
\end{enumerate}
\end{prop}
Note that $\TT:= T_s\SSS$ is independent of
$s\in\SSS$.  Note that $\SSS$ may be parametrized by
either $(\alpha_1,\dots,\alpha_N)$ or 
$(\delta_1,\dots,\delta_N)\in \RR^N\simeq \TT$.
We abuse notation and write
 $s=s(\alpha_1,\dots,\alpha_N)$ or $s=s(\delta_n,\dots,\delta_N)$
to emphasize the dependence.

\Pf
Being a sum of non-negative terms we clearly have 
  $\inf_{\EE} \Om_{0}(f_n,\phi_n,\vec A)\ge 0.$ 
The infimum of zero will be attained  at $(f_n^0,\phi_n^0,\vec A^0)$
if and only if they solve the
following first-order equations in our space $\EE$:
\be\label{firstorder}
f_n(x)\equiv 1, \qquad
 \phi'_n(x)-A_x(x,z_n)\equiv 0 ,\qquad 
 \curl\vec A = H.
\ee
Note that by Lemma~\ref{coulombfinie}  the last
equation is uniquely solved for $\vec A^0\in H^1(D)$, 
with solution as in (\ref{A0}).  By the trace theorem
$A_x(x,z_n)\in H^{1/2}([-L,L])$, and  therefore, $\phi_n$ is
uniquely determined by integration, except for the
 $(N+1)$ constants of integration $\alpha_n$.
  Note that $\eta$ is even in $x$, and hence 
$\int_0^x A^0_x(t,z_n)\, dt$ is an odd function.  In particular,
$\alpha_n$ gives the average value of $\phi_n^0$ and
thus $\alpha_0=0$ is fixed by the definition of the space $\EE$,
leaving $N$ free parameters $(\alpha_1,\dots,\alpha_N)$ to 
parametrize the solution set.  The explicit form for $\Phi^0_{n,n-1}$
is then obtained from (\ref{Phieqn}).
This completes the proof of (a).

Writing the linearized operator around a solution $(f_n^0,\phi_n^0,\vec A^0)$
as a quadratic form,
\bea \nnn
D^2\Om_0 (f_n^0,\phi_n^0,\vec A^0)[u_n,v_n,\vec a] &=&  
  p\sum_{n=0}^N \int_{-L}^L \left\{ 2u_n^2 +
      {1\over\kappa^2}  [u'_n]^2 + {1\over\kappa^2}[v'_n-a_x(x,z_n)]^2
        \right\}\, dx  \\
   \nnn
   && \qquad
   +  {1\over\kappa^2}\iint_D \left| \curl \vec a \right|^2 \, dx\, dz \\
\nnn
   &=& \langle (u_n,v_n, \vec a), (u_n,v_n, \vec a)\rangle -  \\
 \label{quadform2}
   &&\qquad
    p\sum_{n=0}^N \int_{-L}^L 
      \left\{  {2\over\kappa^2}v'_n(x)\, a_x(x,z_n)+v_n^2(x)
          \right\}\, dx
\eea
where $(u_n,v_n,\vec a)\in E$.  

First note that by the second identity in (\ref{quadform2}),
$D^2\Om_0 (f_n^0,\phi_n^0,\vec A^0)$ is of the form
identity plus compact (since the trace embedding of
$\vec a\to \vec a(\cdot, z_n)$ is compact from $H^1(D)$ to
$L^2([-L,L])$.)  Next assume that 
$(u_n,v_n,\vec a)\in\ker D^2\Om_0 (f_n^0,\phi_n^0,\vec A^0)$.
Clearly $u_n\equiv 0$ for all $n=0,\dots,N$.
By Lemma~\ref{control} we must have $\vec a\equiv 0$,
and hence $v_n$ are constant.  Since the definition of the
space $E$ forces $v_0(x)$ with mean zero, we
are left with $N$ free parameters, and 
\beann
  \ker D^2\Om_0 (f_n^0,\phi_n^0,\vec A^0)
  & = &\{ (u_n,v_n,\vec a)\in E: \
           (u_n)_n\equiv 0, \ \vec a \equiv 0, \\
       &&\qquad         (v_0,\dots,v_N)=(0,c_1,\dots,c_N), \
                      c_1,\dots,c_N\in \RR\} \\
        &=& T_s\SSS.
\eeann
\QED

\section{Degenerate perturbation theory.}
\setcounter{thm}{0}

We now perturb away from the degenerate minima of
$\Om_0$, using a variational
Lyapunov--Schmidt procedure.  This method
has been used by Ambrosetti, Coti-Zelati, \& Ekeland \cite{ACE},
 Abrosetti \& Badiale \cite{AB}, Li \& Nirenberg \cite{LN} (and
 many others) in a variety of
situations
 involving heteroclinic solutions of Hamiltonian systems and
 in the semiclassical limit of the nonlinear Schr\"odinger equation.
 
We now proceed according to the usual Lyapunov--Schmidt reduction.
Since $\SSS$ is a hyperplane, $\TT=T_s\SSS$ is independent of
$s\in\SSS$.  Let $W=\TT^\perp$, so any $(f_n,\phi_n,\vec A)\in \EE$ admits the
unique decomposition $(f_n,\phi_n,\vec A)=s+w$ with $s\in\SSS$, $w\in W$,
and any $U:=(u_n,v_n,\vec a)\in E$ decomposes uniquely
as $U=t+w$ with $t\in\TT$, $w\in W$.
We denote the orthogonal projection maps $P: E\to \TT$, $P^\perp:E\to W$
so that $PU=t$, $P^\perp U=w$ whenever $U=t+w$.   
As in the previous
section we interpret the first variation
$\nabla\Omr(f_n,\phi_n,\vec A)$ as an element of $E$ itself, and
 project the equation $\nabla\Omr(f_n,\phi_n,\vec A)=0$
into the two linear subspaces $\TT$ and $W$,
\bea
\label{Teqn}
F_1(r,s,w) & :=  & P\left[ \nabla\Omr (s+w) \right] = 0; \\
\label{Weqn}
F_2(r,s,w) & :=  & P^\perp\left[ \nabla\Omr (s+w) \right] = 0.
\eea
By Proposition~\ref{r=0}~(b),
the second equation can be solved uniquely for
$w=w(r,s)$ in a neighborhood
of $\SSS$ for $r$ small, using the Implicit Function Theorem.
Because our functional $\Omr$ is smooth we can expand 
$w(r,s)$ in powers of $r$.  Note that
$\Omr(s(\alpha_1,\dots,\alpha_N)+w)$ is periodic in 
$(\alpha_1,\dots,\alpha_N)$ so that we may think of $\SSS$ as a (compact)
$N$-torus.  Therefore 
the expansion will be uniform in $s$.  We
obtain the following variant of Lemma~2 of \cite{AB}:
\begin{lem}\label{Lya}
There exist constants $r_0>0$ and $\delta>0$,
depending on $N,L,\kappa,$ and $H$, and a
smooth function
$$  w=w(r,s): \ (-r_0,r_0)\times \SSS\to W\subset E  $$
such that:
\begin{enumerate}
\item[(i)] There exists smooth functions
$w_1,w_2$
such that 
$$  w(r,s)=r w_1(s) + r^2 w_2(r,s)
$$
 for all $|r|<r_0$ and for all $s\in\SSS$;
\item[(ii)]   $P^\perp [\nabla \Omr(s+w(r,s))]=0$.  
\item[(iii)]
Conversely,
if $P^\perp [\nabla \Omr(s+w)]=0$ for some
$r\in(-r_0,r_0)$ and  $w\in W$ with $\|w\|_E < \delta$, then
$w=w(r,s)$.
\item[(iv)]
For any choice of $L_0,\kappa_0,H_0>0$ the constant $r_0$
may be chosen uniformly for all $N\ge 1$, $1\le L \le L_0$,
$1\le \kappa\le \kappa_0$, and $H\ge H_0$.
\end{enumerate}
\end{lem}
Parts (i)--(iii) follow easily from the Implicit
Function Theorem.  The dependences on the various parameters
is more delicate:  we provide the full proof in section~2.5.

\medskip

We define 
$$  \SSS_r := \{s+w(r,s): \ s\in \SSS\}\subset \EE.  $$
$\SSS_r$ is a smooth manifold smoothly diffeomorphic to the
hyperplane $\SSS$.  The important role played by
$\SSS_r$ is that it is a natural constraint for $\Omr$
(see Lemma~4 of \cite{AB}),
and hence the equation (\ref{Teqn}) may be solved variationally:
\begin{lem}\label{natural}     $  $

\noindent
 (a)\
If $(f_n,\phi_n,\vec A)\in \SSS_r$ satisfies
$D(\Omr{}_{|\SSS_r})(f_n,\phi_n,\vec A)=0$, then
$\nabla\Omr(f_n,\phi_n,\vec A)=0$ in $E$.

\noindent
(b)\
There exists $\epsilon_0=\epsilon_0(N,L,\kappa,H)>0$
 such that if 
$(f_n,\phi_n,\vec A)\in \EE$ is a critical point of $\Omr$ with 
\be\label{smallenergy}
   \Omr(f_n,\phi_n,\vec A)\le \epsilon_0,
\ee
 then
$(f_n,\phi_n,\vec A)\in\SSS_r$.

\noindent
(c)\
There exists 
$\tilde r_0= \tilde r_0(N,L,\kappa,H)$ such that for all $0<r<\tilde r_0$,
$\inf_{\EE}\Omr = \inf_{\SSS_r}\Omr$.
\end{lem}
\Pf
The assertion (a) is exactly Lemma~4 of \cite{AB}.  To prove (b),
assume that $(f_n,\phi_n,\vec A)\in \EE$ is a critical point of
$\Omr$ with energy bounded by (\ref{smallenergy}).
 As remarked in the
beginning of the section, $(f_n,\phi_n,\vec A)$ decomposes uniquely as
$(f_n,\phi_n,\vec A)=s+w$,
$s=(f_n^0\equiv 1, \phi_n^0,\vec A^0)\in \SSS$, $w\in W$.  
Now we use the energy bound to estimate the distance between
$(f_n,\phi_n,\vec A)$ and $s$: by Lemma~\ref{control},
\beann
\| \vec A -\vec A_0\|_{H^1} \le
     C_0 \| \curl\vec A - \curl\vec A_0\|_{L^2}  &= & 
       C_0 \|\curl\vec A - H\|_{L^2} \le C_0\epsilon_0 \kappa^2;
         \\
 p\sum_{n=0}^N \int_{-L}^L
    \left[ {1\over \kappa^2} (f'_n)^2 + (1-f_n)^2\right]\, dx
 & \le & 
    p\sum_{n=0}^N  \int_{-L}^L
\left[ {1\over \kappa^2} (f'_n)^2 + (1-f_n^2)^2\right]\, dx \le \epsilon_0; \\
 p\sum_{n=0}^N \int_{-L}^L {1\over \kappa^2}
        (\phi'_n -\phi'_{n0})^2 \, dx
          &=&   p\sum_{n=0}^N \int_{-L}^L {1\over \kappa^2}
        (\phi'_n - A_{x0}(x,z_n))^2 \, dx  \\
            &\le & 2 p\sum_{n=0}^N \int_{-L}^L \left( {1\over \kappa^2}
        (\phi'_n - A_{x}(x,z_n))^2    \right.\\
        &&\quad \left.
          + {1\over \kappa^2}
        (A_x(x,z_n) - A_{x0}(x,z_n))^2 \right)\, dx  \\
       &\le & 2\epsilon_0 + C_1 \|\vec A-\vec A_0\|_{H^1}
         \le C_2\epsilon_0,
\eeann
where $C_1$ comes from applying 
 the Trace theorem  (see Lemma~\ref{trace} in section~5)
in the last line.
Finally, since $[(f_n,\phi_n,\vec A)-s]\in W$,
 each $\phi_n-\phi^0_n$ has mean value zero.
Therefore,  the $H^1([-L,L])$-norm of the difference is
controlled by the difference of the derivatives, as estimated
above, and we may choose $\epsilon_0$ small enough such that
\be\label{close} 
 \mbox{dist}\left( (f_n,\phi_n, \vec A), (f_{n0},\phi_{n0},\vec A_0)
    \right) < c\epsilon_0 <\delta, 
\ee
where $\delta=\delta(N,L,\kappa,H)$ is given by
Lemma~\ref{Lya}.  It then follows by assertion (iii) of Lemma~\ref{Lya}
that if $(f_n,\phi_n,\vec A)$ is a critical point of the Lawrence--Doniach
system with the given energy bound it must lie on $\SSS_r$.
This completes the proof of (b).

To prove (c) we note that  Proposition~\ref{basicenergy} implies that
$\inf_{\EE}\Omr\le 2Np(L+{1\over Hp})r$, and hence we
can choose $\tilde r_0<r_0$ such that 
$2Np(L+{1\over Hp})\tilde r_0\le \epsilon_0$.
\QED

\smallskip

\begin{rem}\label{nodependence}\rm
We note that we cannot make the same statement
about the $N,L,\kappa,H$ dependences of $\tilde r_0$ in
(c) of Proposition~\ref{natural} as we make for $r_0$ in
Lemma~\ref{Lya}.  The uniform bounds on $r_0$ are possible
because of {\it local} estimates on the solution set 
of (\ref{Weqn}) and a continuity argument from $r=0$.
(See section~5.)
We have no such control on the distance $\delta$ from $\SSS$ to
solutions which do {\it not} lie on 
on the manifold $\SSS_r$.  This would entail uniform
(in the parameters) {\it global}
(i.e., non-perturbative) information about the solution set, which 
energy bounds do not provide.
  This leaves open the possibility
of an interval $r\in (\tilde r_0, r_0)$ for which the solutions
on $\SSS_r$ continue to exist and are represented by
a perturbation expansion in $r$, but the absolute minimizer
might not be an element of this family. 
\QED
\end{rem}

In conclusion, we have achieved a complete finite-dimensional
reduction of our problem, for small $r$.  That is to say, when
$0<r<\tilde r_0$ all low energy solutions of the Lawrence--Doniach
system can be found on the $N$-dimensional surface $\SSS_r$.
Moreover, an explicit form for these solutions may be determined
by a simple procedure of Taylor expansion of the equations and
energy in powers of $r$, as is legitimized by Lemma~\ref{Lya}.

\section{Vortex planes.}
\setcounter{thm}{0}

We now apply the theory of the previous section to determine
the minimizer (and other stationary states) of the Lawrence--Doniach
energy, for $r<<1$.  

We summarize our results in the following:

\begin{thm}\label{Thm1}
Assume that \, $\sin(HpL)\neq 0$.
There exists $r_1=r_1(N,L,\kappa,H)$
such that for every $r$ with $0<r<r_1$, the global minimum of
free energy is attained by the vortex plane solutions, given asymptotically
by (\ref{jz})--(\ref{Phi}) below.  Moreover, when $0<r<r_1$
$\Omr$ admits exactly
$2^N$ physically distinct critical points with energy bounded 
as in (\ref{smallenergy}).
\end{thm}

\subsection{Minimizing $\Omr|_{\SSS_r}$}
By Proposition~\ref{natural} we seek
critical points of the finite dimensional variational
problem $\Omr|_{\SSS_r}$.
Using Lemma~\ref{Lya} (i)
we expand a point $(f_n,\phi_n,\vec A)\in \SSS_r$  as
$$  (f_n,\phi_n,\vec A) = s+ w(r,s)=(f_n^0, \phi_n^0, \vec A^0) 
      + r(u_n,v_n,\vec a) + O(r^2),  $$
where $(f_n^0, \phi_n^0, \vec A^0)\in \SSS$ solve the
first-order system (\ref{firstorder}),  and
the error term is uniform over $s\in \SSS$.  
We observe that $\Omr$ has the form,
\beann  \Omr(f_n,\phi_n,\vec A) &=& 
\Omz(f_n,\phi_n,\vec A) + r \Gamma(f_n,\phi_n,\vec A), 
\quad\mbox{with} \\
\Gamma(f_n,\phi_n,\vec A) &=&
      {p\over 2} \sum_{n=1}^N   \int_{-L}^L
                 \left( f_n^2 +f_{n-1}^2 
              -2 f_n f_{n-1}\cos(\Phi_{n,n-1} )\right) \, dx
\eeann
Since $\Omr(s+w(r,s))$
is a smooth function of
$s=s(\delta_1,\dots,\delta_N)$ and $r$, it admits a
Taylor expansion of the form,
\bea  \nnn
\Omr(s+w(r,s)) &=&
    \Omz(s)
       + r\left. {\partial\over\partial r}
             \Omr(s+w(r,s))\right|_{r=0}  + O(r^2)  \\
\nnn
&=&  r\left(
     \Gamma(s) 
     + \nabla\Omz (s)\left[{\partial w\over\partial r}(0,s)\right]\right) +O(r^2)
     = r\, \Gamma(s) +O(r^2)\\
  \nnn
  & = &  rp\sum_{n=1}^N \int_{-L}^L [1-\cos(\delta_n + Hpx)]\, dx  + O(r^2)\\
  \label{est1}
  & = & 
     2p\left( NL- {\sin(HpL)\over Hp}\sum_{n=1}^N\cos \delta_n\right) r +
O(r^2),
  \eea
with remainder term uniform for $s\in \SSS$.

\smallskip

Define $G:\RR\times\RR^N\to \RR$ by:
$$    G(r,\delta_1,\dots,\delta_N):=\Omr(s+w(r,s))/r,  \qquad
     s=s(\delta_1,\dots,\delta_N).  $$
Then $G$ is a smooth function which is periodic in each
coordinate $\delta_n$, and
$$   G(0,\delta_1,\dots,\delta_N)=
    2NpL-{2\sin(HpL)\over H}\sum_{n=1}^N \cos\delta_n. $$ 
Note that when
$$  \sin(HpL)\neq 0 $$
$G(0,s)$ possesses exactly $2^N$ critical points for
$(\delta_1,\dots,\delta_n)\in {\cal K}:=\RR^N/(2\pi{\bf Z})^N$
corresponding to any permutation of 
$$    \delta_n\in \{0,\pi\} \ \mbox{mod $2\pi$.}  $$
It is easy to see that 
 each is a {\it non-degenerate} critical point of $G(0,s)$.
By the Implicit Function Theorem, there exists $\tilde r_1>0$ such
that for all $r\in (0,\tilde r_1)$, and for any
critical point $(\delta^*_1,\dots,\delta^*_N)$ of $G(0,\delta_1,\dots,\delta_N)$
there exists a unique
critical point $(\delta_1(r),\dots,\delta_N(r))$
of $G(r,\delta_1,\dots,\delta_N)$, with
$$   (\delta_1(r),\dots,\delta_N(r))=
         (\delta^*_1,\dots,\delta^*_N) + O(r).  $$
Since $G$ is periodic in each $\delta_n$ we may also conclude 
(via a compactness argument) that
these are the only critical points of $G(r, \delta_1,\dots,\delta_N)$
for small $r$.  Since $\Omr|_{\SSS_r}= rG(r, \delta_1,\dots,\delta_N)$,
by Lemma~\ref{natural}, therefore, for all $r$ with
$0<r<\min\{\tilde r_0,\tilde r_1\}:=r_1$ and whenever $\sin (HpL)\neq 0$,
$\Omr$ admits exactly $2^N$ critical points
(mod~$2\pi$ in each $\phi_n(x)$) with energy bound (\ref{smallenergy}).

\smallskip

The absolute minimizer of $G(0,s)$ is
obtained for
\be\label{vpcond}
  \delta_n^*=  \cases{0, & when ${\sin(HpL)\over Hp}>0$, 
  for every $n=1,\dots,N$;  \cr
\pi, & when ${\sin(HpL)\over Hp}<0$, 
  for every $n=1,\dots,N$, \cr}
  \ee
By the Implicit Function Theorem argument in the previous
paragraph, we conclude that  for all  $r\in (0,r_1)$, 
$\Omr|_{\SSS_r}$ is minimized by a unique
$s_r=s(\delta_1(r),\dots,\delta_N(r))$ with
$$ \delta_n(r) = \delta^* + O(r), \qquad n=1,\dots,N,  $$
where $\delta^*\in\{0,\pi\}$ is
chosen as in (\ref{vpcond}).  Finally, by Lemma~\ref{natural}
we conclude that for all $r\in (0,r_1)$ and $\sin(HpL)\neq 0$
the absolute
minimizer of $\Omr$ in $\EE$ is given 
by $(f_n,\phi_n,\vec A)=s_r+ rw_1(s_r) + O(r^2)$, with 
minimum energy 
$$ \inf_{\EE} \Omr = 2Np \left(
           L-{\sin(HpL\over Hp}   \right)\, r + O(r^2).  $$

\medskip

When $\sin(HpL)=0$ then $\Omr{}_{|\SSS_r}$ is degenerate at order
$r$.  Normally we should go to a higher order in the expansion to determine
the stationary configurations at these exceptional values of
$H=m\pi/Lp$, $m=1,2,3,\dots$, but (as we will see later in the section)
the existence of these degenerate points is explained by the
exchange of stability of two families of solutions when vortices
are nucleated into the sample from the boundary.  We note that
the treatment of the periodic problem in the subsequent paper
\cite{ABB2} will require an expansion of the energy to order
$r^2$ to resolve the degeneracy at any applied field $H$.


\subsection{Expanding the solutions to order $r$}
Because $\Omr$ is a smooth functional we may use
the Implicit Function Theorem and the decomposition
of Lemma~\ref{Lya} to obtain an expansion to arbitrary
order in $r$ of any critical point satisfying the energy
estimate (\ref{smallenergy}).  Here we generate the
expansion to order $r$, to get a better idea of the 
nature of the global minimizers.
Take any such critical point 
of $\Omr|_{\SSS_r}$, with expansion as in (i) of
Lemma~\ref{Lya},
$(f_n,\phi_n,\vec A)= s + rw_1(s) + O(r^2),$
$s=s(\delta_1,\dots,\delta_N)$.  We deduce the equation
satisfied by $w_1(s)$ by implicit differentiation of the equation (\ref{Weqn}):
\bea\nnn
  0 &=& \left.  {d\over dr}
      P^\perp\left(\nabla\Omr(s + w(r,s))\right) \right|_{r=0}
\\
\nnn
      &=&  P^\perp\left(\nabla\Gamma(s) + D^2\Omz(s)
        \left[
           {\partial w\over\partial r}(0,s)
               \right] \right) \\
\label{orderr}
   &=& P^\perp\left(\nabla\Gamma(s) + D^2\Omz(s)
        \left[
           w_1(s)
               \right] \right).
\eea
Since $D^2\Om_0(s)$ is an invertible map from $W\to W$
this formula uniquely determines $w_1(s)$.  Now we represent
$w_1(s)=(u_{n,1},v_{n,1},a_{x,1}, a_{z,1})$, in other words
\beann
   f_n=1+ru_{n,1} + O(r^2),  &\quad & 
     \phi_n= \alpha_n + nHpx + v_{n,1}+O(r^2),  \\
           A_x= Hz + ra_{x,1} + O(r^2), &\quad &
            A_z=ra_{z,1}+O(r^2), 
\eeann
and denote $b(x,z)=\curl\vec a_1=\partial_z a_{x,1}-\partial_x a_{z,1}$.
In terms of these coordinates (\ref{orderr}) takes the following
simple form:
$$  
-{1\over \kappa^2} u''_{n,1} + 2u_{n,1}  = {1\over 2}\cases{
          \cos(\delta_n+Hpx)
            +\cos(\delta_{n+1} + Hpx)-2,
                & $1\le n\le N-1$,  \cr
        \cos(\delta_N + Hpx)- 1, 
           & $n=N$, \cr
         \cos(\delta_1 + Hpx)- 1, 
           & $n=0$, \cr }
$$
with $u'_{n,1}(\pm L)=0$;
\beann
&&  {1\over \kappa^2}{d\over dx}\left( v'_{n,1} - a_{x,1}(x,z_n)\right) = 
   \\
 &&\qquad      {1\over 2}\cases{
          \sin(\delta_n + Hpx)- \sin(\delta_{n+1}+Hpx) - I_n, 
                & $1\le n\le N-1$, \cr
        \sin(\delta_N + Hpx)-I_N, & $n=N$, \cr
       -\sin(\delta_1+Hpx)-I_0, & $n=0$, \cr }
\eeann
with boundary condition $v'_{n,1}(\pm L) - a_{x,1}(\pm L, z_n)=0$
and 
$$  I_n={1\over 2L}
   \cases{   \int_{-L}^L [\sin(\delta_n + Hpx)- \sin(\delta_{n+1}+Hpx)]
           \, dx, & $1\le n\le N-1$, \cr
       \int_{-L}^L \sin(\delta_N + Hpx)\, dx, & $n=N$, \cr
       - \int_{-L}^L\sin(\delta_1+Hpx)\, dx, & $n=0$; \cr }
$$
and
\be
\label{proj3} 
\vec a_1(x,z)= \left({\partial\xi\over\partial z}, 
       -{\partial\xi\over\partial x}\right),
          \qquad \Delta\xi=b(x,z), 
              \quad \left.\xi \right|_{\partial B}=0,
\ee
where  $b(x,z)=b^{(n)}(x)$ for $z_{n-1}<z<z_n$, with
\be
\label{proj4}
    {\partial b^{(n)}\over \partial x} =
          {p\kappa^2\over 2} \left[
      \sin(\delta_n + Hpx)-{1\over 2L}\int_{-L}^L \sin(\delta_n+Hpx)\, dx\right],
            \qquad
      b^{(n)}(\pm L)=0,
\ee
for $n=1,\dots,N$.

\medskip

Now assume that $\sin(HpL)\neq 0$ and
consider the absolute minimizers, $\delta_n\equiv\delta_*$
with $\delta_*\in \{0,\pi\}$ chosen according to (\ref{vpcond}).
 Except for an edge effect at the top and bottom layers
($n=0,N$) the gauge-invariant quantities are independent
of $n$ at order $r$.  In particular,  the magnetic field $h$,
and Josephson current density $j_z^{(n)}$ are (at order $r$)
completely independent of $z$, $n$:
\bea
\label{jz}
&&j_z^{(n)}:=j_c f_n f_{n-1}\sin\Phi_{n,n-1}= 
r {\kappa^2 p\over 2} \sin(\delta+Hpx) + O(r^2), \\
\label{h}
&&h(x,z)=H+r{\kappa^2\over 2 H}\left(\cos(\delta+HpL)-
\cos(\delta+Hpx)\right) + O(r^2)
\eea
where $j_c=r {\kappa^2 p\over 2}$ is the critical Josephson current.
The in-plane supercurrent vanishes at order $r$ for all interior planes,
$$     j_x^{(n)}=O(r^2), \qquad  n=1,\dots,N-1,  $$
but the top and bottom of the sample carry current at order $r$:
\bea
\label{top}
&&j_x^{(N)}=0+r{\kappa^2\over 2Hp}\left(
\cos(\delta+HpL)-\cos(\delta+Hpx)\right) + O(r^2)\\
\label{bottom}
&&j_x^{(0)}=0+
r{\kappa^2\over 2Hp}\left(
\cos(\delta+Hpx)-\cos(\delta+HpL)\right) + O(r^2)
\eea
Similarly, the order $r$ expansion of modulus of the order parameter is 
$n$-independent on the interior planes, but is modified
at the top and bottom at order $r$:
\be
\label{fn}
f_n(x)=1+r\left(-{1\over 2}+ A\cosh({\sqrt 2\over \kappa}x)+
B\cos(\delta+Hpx)\right) + O(r^2),  \qquad n=1,\dots,N-1,
\ee
\be
f_0=f_N=1+{r\over 2}\left(-{1\over 2}+ A\cosh({\sqrt 2\over \kappa}x)+
B\cos(\delta+Hpx)\right) + O(r^2)
\ee
where 
$$   A={\kappa^3Hp\over\sqrt
2(H^2p^2+2\kappa^2)}{\sin(\delta+HpL)\over\sinh({\sqrt 2 L\over\kappa})},
\quad B={\kappa^2\over H^2p^2+2\kappa^2}.$$
Expansion to order $r^2$ in the solutions will show these quantities to
be independent of $z$ (or $n$) except for the top and bottom {\it two}
planes and the top and bottom gaps.  

We note that the gauge-invariant phase difference,
\be \label{Phi}
\Phi_{n,n-1}(x)= \delta+Hpx + r {\kappa^2\over 2 H^2}
\left(xHp\cos(\delta+HpL)-\sin(\delta+Hpx)+c_n\right) + O(r^2), 
\ee
contains constants $c_n$ which can only be determined by higher
order expansion in $r$.  Indeed, they represent the order $r$
correction to the choice of $\delta_n$ when minimizing the
finite dimensional problem $\Omr|_{\SSS_r}$. 

We call this configuration
{\it vortex planes}-- see Figure 1.  
Unlike the Ginzburg--Landau case, the order
parameter does not vanish at the ``core'' where the local field
attains its maximal value and around which supercurrents circulate.
The ``vortices'' are then the planes $\{x=const.\}$ over which
$h(x)$ attains its relative maxima.  These planes are nodes for
the current, and flux per plane per cycle of the Josephson current
is, to order $r$, the usual flux quantum ($2\pi$ in our units.)

\begin{rem}\label{surface}\rm \  (a)\
The distinction between the energies of the various
lattice geometries at the lowest order
term $O(r)$ in the expansion (\ref{est1}) is a {\it surface term},
in the sense that it scales like the length of the lateral
edges of the sample $2Np$, as opposed to the free energy
itself which scales like the cross-sectional area of the
bulk $2LNp$.  In other words, the vortex plane configuration
is preferred for the effect it has on the surface currents on
the left and right edges of the sample, with no regard to
energy minimization in the interior.  As we will see in 
the periodic case in \cite{ABB2} the effect on the bulk will
be observed in an order $r^2$ term.  For all other parameters
fixed, eventually we can take $r>0$ small enough so that
the $O(r)$ surface term dominates. But in a real superconductor
the value of $r$ is given, and hence when the width $L$ is
increased eventually the bulk $O(r^2)$ term will compete with the
surface $O(r)$ term.  This suggests that the interval of validity
of the $r$ expansion may not be uniform, but rather deteriorate
with increasing sample width $L$.

\noindent
(b)\
Another indication that the radius of convergence of the expansion 
depends inversely on the width $L$ is the presence of a linearly
growing factor (secular term) in the order $r$ term of the expansion
(\ref{Phi}) of $\Phi_{n,n-1}(x).$ 
\QED
\end{rem}

\subsection{Vortex plane nucleation}
The very precise desciption of minimizers for $r\simeq 0$
allows us to identify the transitions which the
sample undergoes as new vortices are nucleated in
an increasing applied field $H$.
When $0\le H<\pi/pL$ we observe the Meissner state for this problem:
$h$ attains its maximum at the boundary $x=\pm L$, and decreases to a line of
minima at $x=0$.  For our solutions the lower critical field
$H_{C1}=\pi/pL$:  at this point the finite dimensional minimization
problem $G(0,s)=0$ degenerates, and when $H$ is increased slightly
the minimizing configuration has phases
 $\delta_n=\pi$ for all $n=1,\dots,N$.  Note that
at this critical value of $H$ a new node appears in the Josephson current 
(\ref{jz}) at each endpoint.
 The switch in $\delta_n$ amounts to
a flip in sign for the Josephson currents and for the variation
of $h$ from $H$, so the newly nucleated nodes correspond to
{\it minima} of $h$, with the former minimum at $x=0$ becoming
a local {\it maximum}.  The same phenomenon will occur every
time $H$ crosses a value $k\pi\over pL$, $k=2,3,\dots$
Each time, two new nodes for the Josephson current will be
nucleated, and the change of $\delta_n$ from zero to $\pi$ (or
vice-versa) will exchange minima and maxima of $h$ in the
interior, resulting in the creation of exactly {\it one} new
vortex plane.  Since the minimum energy is given asymptotically
by
$$  \epsilon(H):=\min \Omr = 
      Npr\left(2L-\left|{2\sin(HpL)\over Hp}\right| \right)
           + O(r^2),  $$
 the magnetization $M(H):=\partial\epsilon/\partial H$ is discontinuous
 at each nucleation, indicating a first-order phase transition.

\section{The validity of the expansion}
\setcounter{thm}{0}

In this section we prove Lemma~\ref{Lya} which justifies the
finite dimensional reduction in a neighborhood of $\SSS$.
While parts (i)--(iii) follow easily from the
Implicit Function Theorem, the real interest is in the
dependence of the interval of validity $|r|<r_0$ on the
many parameters in the model, especially the sample
dimensions $L$ and $M=Np$, $\kappa$, and the applied field $H$.
Here we
give a lower bound for the ``radius of convergence'' in
$r$ as a function of these parameters.  We show that the small $r$
approximation is essentially independent of the number of planes $N$, but
that its validity can be expected to decrease with increasing width
$L$ or increasing $\kappa$.  On the other hand, increasing
the external field $H$ enhances the approximations somewhat.
This lower bound is consistent with experiments on the
high-$T_c$ materials, where the vortex planes have not
been observed, perhaps because the large values
of $\kappa$ and the large size of typical samples
(measured in terms of $\lambda_{ab}$)  reduce the radius of validity
of the expansion below the appropriate value of $r$ for
such materials.   For other types of layered superconductors
with smaller values of $\kappa$  we are more likely
to see configurations similar to the solutions produced
in our $r\to 0$ limit.  

\smallskip

Throughout this section, we concentrate on the parameters
$N,L,\kappa,H$, and assume $0<p<1$, $\kappa\ge 1$, and $L\ge 1$.
We also define a norm on the space $W=[T\SSS]^\perp$:  let
$w=(u_n,v_n,\vec a)\in W$.  Then we denote
$$
\|w\|^2 = \|[(u_n,v_n,\vec a)\|^2 =
    p\sum_{n=0}^N   \int_{-L}^L
    \left( (u'_n)^2 +
     u_n^2 + (v'_n)^2 + v_n^2 \right)\, dx
          + \iint_D\left[ |\nabla\vec a|^2 + |\vec a|^2\right]\, dx\, dz.
$$

\medskip

Our approach is to calculate {\it a priori} estimates on
the solutions to the equation (\ref{Weqn}) to determine
when that equation can degenerate.  As a first step,
we must recognize the system of differential equations
satisfied by the solutions to (\ref{Weqn}).  
For any $U=(f_n, \phi_n,\vec A)$ we write
$U= s + w(r,s)$, where
$s=s(\delta_1,\dots,\delta_N)\in \SSS$, and $w(r,s)=(u_n,v_n,\vec a)\in W$.
We recall the effect that this decomposition has on
some familiar quantities:
\beann
&& f_n=1+u_n,  \qquad
 V_n(x)=\phi'_n-A_x(x,z_n)=v'_n-a_x(x,z_n),    \\  \\
&&   \vec  a(x,z)=\left(\partial_z \xi, -\partial_x \xi\right),
               \qquad  \Delta\xi = b(x,z), \quad \xi|_{\partial D}=0,  
 \qquad h(x,z)=H + b(x,z), \\ 
\\ &&    \Phi_{n,n-1}=\phi_n-\phi_{n-1}-\int_{z_{n-1}}^{z_n} A_z (x,t)\, dt
     = Hpx + \delta_n + \varphi_{n,n-1},  \\
&&  \qquad\mbox{where} \quad  \varphi_{n,n-1}(x) = v_n(x)- v_{n-1}(x) - 
     \int_{z_{n-1}}^{z_n} a_z (x,t)\, dt. 
\eeann
The equations for $(u_n,v_n,\vec a)$ will differ from the
unconstrained Euler--Lagrange equations derived in section~2
because of Lagrange multipliers created by projection into the subspace
$W$.  
     
\smallskip

To say $w(r,s)$ solves
$P^\perp D\Omr (s+w(r,s))=0$ is equivalent to saying that
 $\Omr(s+w(r,s))$
is stationary with respect to variations in $W$.
Taking the variation with respect to just one 
of the $u_n$ yields the exact
same equations (\ref{feqn}) as the unconstrained case, which we rewrite in
terms of $u_n$:
\bea
\label{proj1}  && 
-{1\over \kappa^2} u''_n + (u_n+1)(u_n+2)u_n 
    + {1\over\kappa^2}(v'_n-a_x(x,z_n))^2 u_n = 
       -{1\over\kappa^2}(v'_n-a_x(x,z_n))^2 \\
    \nnn && \qquad\qquad
    +\cases{
          {r\over 2}\left(
       (1+u_{n-1})\cos\Phi_{n,n-1}
      +(1+u_{n+1})\cos\Phi_{n+1,n}-2(1+u_n)\right), & $n \neq 0,N$; \cr
  &\cr
     {r\over 2}\left((1+u_1)\cos\Phi_{1,0}-(1+u_0)\right), & $n=0$; \cr
  &\cr
      {r\over 2}\left(
       (1+u_{N-1})\cos\Phi_{N,N-1}
         -(1+u_N)\right),  & $n=N$,\cr} 
\eea
with boundary conditions  $u'_n(\pm L) = 0$.

For $v_n$, $n=0,\dots,N$ the integral constraint gives rise
to the usual Lagrange multiplier, without which  the equations 
generally could not be integrated:
\bea
\label{proj2} &&  
{1\over \kappa^2}{d\over dx}\left( f_n^2(v'_n - a_x(x,z_n))\right) = \\
\nnn && \qquad\qquad
   \cases{
       {r\over 2}[f_n  f_{n-1}\sin\Phi_{n,n-1}
       -f_{n+1}f_n \sin\Phi_{n+1,n}]-C_n, & $n=1,\dots,N-1$;\cr &\cr
     {r\over 2} f_{N} f_{N-1} \sin\Phi_{N,N-1}-C_N,
          & $n=N$, \cr & \cr
     -{r\over 2} f_0 f_1 \sin\Phi_{1,0}+C_0,
          & $n=0$, } 
\eea
where
$$
C_n = \cases{  {1\over 2L}\int_{-L}^L
      {r\over 2}[f_n  f_{n-1}\sin\Phi_{n,n-1}
       -f_{n+1}f_n \sin\Phi_{n+1,n}]\, dx, & $n=1,\dots,N-1$; \cr
       {1\over 2L}\int_{-L}^L
       {r\over 2} f_{N} f_{N-1} \sin\Phi_{N,N-1}\, dx,
          & $n=N$, \cr
       {1\over 2L}\int_{-L}^L
       {r\over 2} f_0 f_1 \sin\Phi_{1,0}\, dx,
          & $n=0$.   }
$$
The equations for $v_n$ also include a no-flux boundary condition,
$$  \left. f_n^2 (v'_n - a_x(\cdot, z_n))\right|_{x=\pm L}=0,
     n=0,\dots,N.  $$

Finally, we derive the Euler--Lagrange
equations for $\vec a$.  As in our previous calculations
it will be easier to deal with the associated magnetic
field $b(x,z)=\curl \vec a(x,z)$.  
\begin{lem}\label{bfield}
If $(f_n,\phi_n,\vec A)=s+w(r,s)$ satisfies (\ref{Weqn}),
then $b(x,z)= h(x,z)-H$ satisfies:
\be
\label{proj4a} 
  b(x,z)= b^{(n)}(x), \quad z_{n-1}<z<z_n, \quad n=1,\dots,N, 
\ee
\be\label{proj4b}
   b^{(n)}(x)-b^{(n+1)}(x) = p(v'_n - a_x(x,z_n)), \quad n=1,\dots,N-1, 
\ee
\be\label{proj4c}
     b^{(N)}(x)
            =   p \left( (v'_N-a_x(x,z_N)) f_N^2\right), \qquad
     b^{(1)}(x)
            =  - p\left( (v'_0-a_x(x,z_0)) f_0^2\right) 
\ee
\bea\label{proj4d}
   {d b^{(n)}\over d x } &=&
          {rp\kappa^2\over 2} f_n f_{n-1}\sin\Phi_{n,n-1}(x)
              - D_n, \quad
b^{(n)}(\pm L)=0,  \\
\nnn  \mbox{where}\  D_n &=& {1\over 2L}\int_{-L}^L 
        {rp\kappa^2\over 2} f_n f_{n-1}\sin\Phi_{n,n-1}(x)\, dx, 
\eea
\end{lem}
\Pf
First we derive the equation for $\vec a$ in the weak form.
Any admissible variation for $\vec a$ in $E$ 
can be represented in the form
$\tilde a=(\partial_z\eta, -\partial_x \eta)$
with $\eta\in H^2\cap H^1_0(D)$.  (Indeed, $\eta$ solves
$\Delta\eta=\curl \tilde a$ in with $\eta|_{\partial D}=0$:
see Lemma~\ref{control}.)
Taking the first variation of energy in the direction
of $\tilde a$ we obtain:
\bea\label{weaka}
\nnn
0 &=& -p\sum_{n=0}^N \int_{-L}^L
      {1\over\kappa^2}(v'_n-a_x(x,z_n))f_n^2 \partial_z\eta(x,z_n) \, dx \\
&&\quad  + \sum_{n=1}^N \int_{-L}^L\int_{z_{n-1}}^{z_n}
      \left\{
        {1\over\kappa^2} b(x,z) \Delta\eta +
           {rp\over 2} f_n f_{n-1}\sin\Phi_{n,n-1}\, \partial_x\eta
             \right\}\, dz\, dx  , 
\eea
for all $\eta\in H^2\cap H^1_0(D)$. 

\smallskip

Next, we observe that there exist functions $\hat b^n(x)\in H^1(-L,L)$
which solve the system (\ref{proj4b})--(\ref{proj4d}).
Indeed, we define $\hat b^n$ by integration,
$$  \hat b^n(x) = 
          \int_{-L}^x \left( 
   {rp\kappa^2\over 2} f_n(x) f_{n-1}(x)\sin\Phi_{n,n-1}(x)
     - D_n\right)\, dx, \quad n=1,\dots,N.  $$
Then $\hat b^n(x)\in H^1([-L,L])$ and satisfies (\ref{proj4d}).
For $n=1,\dots,N-1$ we have:
\beann
 {d\over dx}
        \left( \hat b^{n+1}(x)-\hat b^n(x) \right)
            &=&  {rp\kappa^2\over 2}
                \left(  f_{n+1} f_n\sin\Phi_{n+1,n}
                   - f_n f_{n-1} \sin\Phi_{n,n-1} \right)  - D_{n+1}
                     +D_n  \\
         &=& - p{d\over dx} \left( (v'_n-a_x(x,z_n)) f_n^2\right),
\eeann
by (\ref{proj2}).
For the top and bottom gaps ($n=1,N$) we obtain:
\beann
  {d\over dx} \hat b^N(x)
            &=& {rp\kappa^2\over 2}
                   f_N f_{N-1} \sin\Phi_{N,N-1}  
                      -D_N \\
               & =&  p{d\over dx} \left( (v'_N-a_x(x,z_N)) f_N^2\right),\\
           {d\over dx} \hat b^0(x)
            &=&  {rp\kappa^2\over 2}
                   f_0 f_1 \sin\Phi_{1,0}  
                      -D_0 \\
               & =& - p{d\over dx} \left( (v'_0-a_x(x,z_0)) f_0^2\right)
\eeann
In each case we integrate (and use the boundary
conditions $\hat b^n(-L)=0=v_n'(-L)-a_x(-L,z_n)$) to verify
that $\hat b^n$ satisfies conditions (\ref{proj4b}) and
(\ref{proj4c}).

Next we define $\hat b(x,z)=\hat b^n(x)$ for $z_{n-1}<z<z_n$, and
show that $\hat b$ also solves (\ref{weaka}). 
By an integration by parts in each summand,
\beann
&&  \sum_{n=1}^N \int_{-L}^L\int_{z_{n-1}}^{z_n}
        {1\over\kappa^2}\hat b(x,z) \Delta\eta \, dz\, dx  
      =  \sum_{n=1}^N \int_{-L}^L\int_{z_{n-1}}^{z_n}
        {1\over\kappa^2}\hat b^n(x) \Delta\eta \, dz\, dx  \\
&& \qquad\qquad =   {1\over\kappa^2}\sum_{n=1}^N \left\{
          \int_{-L}^L  \hat b^n(x) 
                \left( \partial_z\eta (x,z_n) - \partial_z\eta(x,z_{n-1})
                     \right)\, dx  \right.  \\
& &\qquad\qquad \qquad\left.  -  \int_{-L}^L\int_{z_{n-1}}^{z_n}
          {d \hat b^n\over dx}\, \partial_x\eta(x,z) \, dz\, dx \right\} \\
&& \qquad\qquad = {1\over\kappa^2}
      \int_{-L}^L \left[ \hat b^N(x)\partial_z\eta(x,z_N) -
                     \hat b^1(x)\partial_z\eta(x,0) \right]\, dx  \\
&&\qquad\qquad \qquad + {1\over\kappa^2}\sum_{n=1}^{N-1} 
           \int_{-L}^L  \left( \hat b^n(x) -\hat b^{n+1}(x) \right)
                                        \partial_z\eta(x,z_n) \, dx    \\
&& \qquad \qquad\qquad -\sum_{n=1}^N
          \int_{-L}^L\int_{z_{n-1}}^{z_n}
             \left({rp\over 2}  f_n f_{n-1}\sin\Phi_{n,n-1}(x) -D_n
                  \right)\partial_x\eta \, dz\, dx  \\
&& \qquad\qquad = p\sum_{n=0}^N \int_{-L}^L
       {1\over\kappa^2}(v'_n - a_x(x,z_n))f_n^2 \partial_z \eta(x,z_n)\, dx \\
&& \qquad\qquad\qquad -\sum_{n=1}^N
          \int_{-L}^L\int_{z_{n-1}}^{z_n}
             {rp\over 2}  f_n f_{n-1}\sin\Phi_{n,n-1}(x)
                  \partial_x\eta \, dz\, dx .
\eeann
for all $\eta\in H^2\cap H^1_0(D)$.
That is, $\hat b$ also solves the equation (\ref{weaka}).
Now consider the function $g(x,z)$ defined by
$g(x,z) = b(x,z) - \hat b(x,z).$  Then $g\in L^2(D)$ and
solves
$$  \iint_D g(x,z)\Delta\eta \, dx\, dz =0  $$
for all $\eta\in H^2\cap H^1_0(D)$.  Taking
$\eta$ to be the solution to the Dirichlet problem,
$\Delta\eta=g$, $\eta|_{\partial D}=0$ we arrive at the
desired conclusion $g\equiv 0$, and hence
$b(x,z)=\hat b(x,z)$ and 
the lemma is established.
\QED

\bigskip

As before, we observe that the gauge-invariant phase
difference satisfies a useful equation (see (\ref{Phieqn})),
which derives from Stokes' Theorem:
\be\label{varphieqn}   {d\over dx}\varphi_{n,n-1}(x)
   = \left( v'_n-a_x(x,z_n) \right) - 
                \left( v'_{n-1}-a_x(x,z_{n-1}) \right)
                   + p b^{(n)}(x).  
\ee

\bigskip

We  now begin the proof of Lemma~\ref{Lya}.

\noindent {\bf Step 1:}\
Getting started.
By the implicit function theorem, applied to the equation
(\ref{Weqn}), for every $s\in \SSS$, 
and for every fixed choice of parameters
$N,L,\kappa,H$, there exist constants
$\rho>0$, $r_1>0$ such that (\ref{Weqn}) admits a unique
solution $w=w(r,s)$ with $|r|<r_1$ and $\|w\|<\rho$.
Since $\Omr(s+w(r,s))$ is periodic in $s$ we may compactify
the problem by treating $\SSS$ as an $N$-torus, and 
hence the constants $r_1,\rho$ may be chosen independently
of $s\in \SSS$.  
Statements (i)--(iii) of Lemma~\ref{Lya} then follow from standard arguments
involving the Implicit Function Theorem and the regularity
of the functional $\Omr$
(see \cite{AB}.)  What remains to prove is the dependence of
the value $r_1$ on the various parameters.

We first observe that Lemma~\ref{Lya}~(i) implies that
the solutions $w(r,s)$ satisfy a uniform (in $s$) estimate
$\|w(r,s)\|\le C_1 |r|$, with constant $C_1$ possibly
depending on all parameters, $N,L,\kappa, H, p$. 
From this preliminary argument we may, by
reducing $r_1>0$ if necessary, assume that
\be\label{fnotzero}
   \frac12 \le f_n(x)\le \frac32, \quad n=0,\dots,N  \ee
for $|r|<r_1$.

\medskip

\noindent {\bf Step 2:}\ Uniform estimates in $n$.
We now use the Euler--Lagrange equations 
(\ref{proj1})--(\ref{proj4d}) derived above
to estimate the gauge-invariant quantities,
in the supremum norm and individually in each plane or gap.

First we make rough bounds, then iterate to improve the dependence
on parameters.  Since we can estimate the constants $|C_n|\le r$ and
$|D_n|\le r\kappa^2 p/2$, we have
\bea \label{simipleest1}
  \|b^{(n)}\|_\infty &\le& c r\kappa^2 p L, \\
  \label{simpleest2}
  \| v'_n-a_x(\cdot,z_n)\|_\infty &\le&
    cr\kappa^2 L.
\eea
(Here and in the following $c$ will denote a pure constant,
independent of all other parameters.)

Now consider the equations for $u_n$.  Suppose first that
$\max_{x,k} u_k(x) = u_n(x_0)\ge 0$.  If $n\in \{1,\dots,N-1\}$,
we observe that the extremal property of $u_n(x_0)$ implies
that
$$  u_{n-1}(x_0)\cos\Phi_{n,n-1} + u_{n+1}(x_0)\cos\Phi_{n+1,n} -2 u_n(x_0)
   \le 0, $$
while if $n=0$ or $n=N$ the corresponding terms
$$   u_{N-1}(x_0)\cos\Phi_{N,N-1} - u_N(x_0) \le 0, \ \mbox{or}\
  u_1(x_0)\cos\Phi_{1,0} - u_0(x_0) \le 0.  $$
Consequently, at the maximum point the right-hand side of
(\ref{proj1}) is non-positive.  By the maximum principle
we obtain $u_n(x_0)<0$.  (Note this argument holds in
general, for any value of the parameters, for any
solution of the Lawrence--Doniach system.)

Next assume that
$\min_{x,k} u_k(x) = u_n(x_0)< 0$.   By the same arguments
as above at the point $x_0$, the extremal
$u_n$ satisfies an equation of the form
$$   -{1\over\kappa^2}u''_n(x_0) + c(x_0)u_n(x_0)
\ge -  {1\over\kappa^2} (v'_n -a_x(x_0,z_n))^2 - 2r,  $$
with $c(x)\ge \frac12$.  In conclusion we may bound
$u_n$ via
\be\label{ubound}  0> u_n(x) \ge - {2\over \kappa^2}
    \sup_x (v'_n -a_x(x_0,z_n))^2 - 4r.  
\ee
Applying  the simple estimate (\ref{simpleest2}) we obtain
a preliminary bound on $u_n$,
$$  0> u_n(x) > -cr\left( r\kappa^2 L^2 + 1 \right).  $$
We use the calculus inequality
$$   \sup_{x\in [-L,L]} |f'(x)|
   \le {2\over\delta} \sup_{x\in [-L,L]} |f(x)| 
    + {\delta\over 2}\sup_{x\in [-L,L]} |f''(x)|,
\quad \mbox{for any $0<\delta\le L$,}  $$
with $\delta=1$, and obtain
$$  \|u'_n\|_\infty \le cr\left( r\kappa^2 L^2 + 1 \right).  $$

\smallskip

To improve the above estimates we need  to take into account
the fact that the integrands and right-hand sides contain
terms which oscillate rapidly when $H$ is large.  In particular,
$\Phi_{n,n-1}= Hpx + \varphi_{n,n-1}$, and we expect that
terms containing the sine or cosine of $\Phi_{n,n-1}$ 
tend to zero in the weak sense as $H\to\infty$.
Since these oscillatory terms are sources in the equations
for the fields and currents, we should obtain stronger
convergence to zero as $H\to\infty$, and hence sharper
estimates for large $H$ than we obtained in bluntly
measuring these terms in the  supremum norm.
With this in mind, we require one further ingredient: a
bound for $\varphi'_{n,n-1}$ using equation (\ref{varphieqn}),
\be\label{phi'bound}
\sup_x |\varphi'_{n,n-1}(x)| \le cr\kappa^2 L (1+p^2).
 \ee

We now integrate by parts in $C_n$, $D_n$, and in the integrals
representing $b^{(n)}(x)$ and $v'_n-a_x(x,z_n)$: 
\bea \label{bestb}
\sup_x |b^{(n)}(x)| &\le &
   c{r\kappa^2 p\over Hp} \left(1 + L \|u'_n\|_\infty\right)
       + L \|\varphi'_{n,n-1}\|_\infty \\
 \nnn      &\le & 
    c{r\kappa^2 p\over Hp} ( 1+ rL^2\kappa^2)(1+ rL), 
\eea
and similarly
\be\label{bestj}  \sup_x |v'_n(x)-a_x(x,z_n)|
 \le c{r\kappa^2 \over Hp} ( 1+ rL^2\kappa^2)(1+ rL).
\ee
Substituting into (\ref{ubound}) we obtain a more refined estimate
for $u_n$,
\be\label{bestu}  0> u_n(x) > -  C_u r\left(   1+
     r{\kappa^2 \over H^2 p^2}(1+r\kappa^2 L^2)^2 (1+rL)^2
     \right),
\ee  
with constant $C_u$ independent of $N,\kappa,L$.

\bigskip

\noindent {\bf Step 3:}\ A lower bound on the linearization
at $r=0$.
We define
$$  \lambda=\lambda(N,L,\kappa,p):=
   \inf\{ Q_0(\Psi,\Psi): \  \Psi=(\mu_n,\nu_n,\vec\alpha)
      \in W, \  \|(\mu_n,\nu_n,\vec\alpha)\|=1\},
$$
with
\beann   Q_0(\Psi,\Psi)&=& D^2\Om_{r=0}(f_n^0,\phi_n^0,\vec A^0)
    [\Psi,\Psi] \\
   &=&  p\sum_{n=0}^N \int_{-L}^L \left(
       {1\over\kappa^2} (\mu'_n)^2 + 2\mu_n^2
          + {1\over\kappa^2} \left( \nu'_n -\alpha_x(x,z_n)^2\right)
          \right)\, dx \\
&&\qquad + {1\over\kappa^2} \int_{-L}^L\int_0^{Np}
                 (\curl\alpha)^2\, dz\, dx.
\eeann
Recall that, as an operator on the whole tangent space $E$,
$D^2\Omz(s)$ has an $N$-dimensional kernel for any
$s\in\SSS$.  Here we consider its restriction to the orthogonal
subpace $W=T\SSS^\perp=\ker(D^2\Omz(s))^\perp$.
From the proof of Proposition~\ref{r=0}~(b)
we may conclude that $\lambda>0$ for
any choice of $N,L,\kappa,p$.  Here we will obtain the more
precise information on its dependence on these parameters.

Note first that (trivially, since $\kappa\ge 1$ by hypothesis,) 
$$   p\sum_{n=0}^N \int_{-L}^L \left( {1\over\kappa^2}(\mu'_n)^2 
     + \mu_n^2\right) \, dx  \ge {1\over\kappa^2} 
    p\sum_{n=0}^N \|\mu_n\|_{H^1}^2.  $$

We require the following lemma to estimate the terms
including the vector potential $\vec a$ appearing in
the linearizations:
\begin{lem}\label{trace}
For any $\vec a\in H^1(D)$ with $\vec a\cdot \hat n=0$
on $\partial D$, 
\be\label{tracex}   p\sum_{n=0}^N \|a_x(\cdot,z_n)\|_{L^2(D)}^2 
   \le {(p+1)(N+1)\over N} \| a_x \|_{H^1(D)}^2, 
\ee
\be\label{tracez}  \sum_{n=1}^N \|\int_{z_{n-1}}^{z_n} a_z(z)\, dz\|_{L^2(D)}^2
    \le p \iint_D a_z^2(x,z)\, dx\, dz.  
\ee
\end{lem}

\Pf
Let $\chi_n^+(z) = (z-z_{n-1})/p$, $n=1,\dots,N$.
Then we apply the divergence theorem to the vector
field $(a_x^2(x,z)\chi_n^+,0)$ in the strip $[-L,L]\times (z_{n-1},z_n)$
to obtain
\bea\label{trace1}
\int_{-L}^L (a_x^2(x,z_n))^2\, dx
   &=&  \int_{-L}^L \int_{z_{n-1}}^{z_n}
            {\partial\over\partial z}\left(
                  \chi_n^+(z) a_x^2(x,z) \right) \, dz\, dx \\
\nnn    &\le & \left[ 1+{1\over p}\right]\int_{-L}^L \int_{z_{n-1}}^{z_n}
         \left( a_x^2 
               + \left| {\partial a_x\over\partial z}\right|^2
                \right)\, dz\, dx
\eea
for $n=1,\dots,N$.  Similarly, using $\chi_n^-=(z_{n+1}-z)/p$
in the strip $[-L,L]\times (z_n,z_{n+1})$
we obtain
\be\label{trace2}
\int_{-L}^L (a_x(x,z_n))^2\, dx \le
   \left[1 +\frac{1}{p} \right]
      \int_{-L}^L  \int_{z_{n}}^{z_{n+1}}
         \left( a_x^2 
               + \left| {\partial a_x\over\partial z}\right|^2
                \right)\, dz\, dx,
\ee
for $n=0,\dots,N-1$.

Clearly there exists an index $n_0$ such that
$$  \int_{-L}^L  \int_{z_{n_0-1}}^{z_{n_0}}
         \left( a_x^2 
               + \left| {\partial a_x\over\partial z}\right|^2
                \right)\, dz\, dx 
                \le {1\over N} \iint_D \left( a_x^2 
               + \left| {\partial a_x\over\partial z}\right|^2
                \right)\, dz\, dx.
$$
For $0\le n\le n_0-1$ we use estimate (\ref{trace2})
for $\| a_x\|_2^2$, and for $n_0\le n\le N$ we use
(\ref{trace1}).  In doing so we require the $H^1$-norm of
$a_x$ in each interval exactly once, except for the
strip $z_{n_0-1}<z<z_{n_0}$ which appears {\it twice.}
In this way we arrive at the desired estimate (\ref{tracex}).

The estimate (\ref{tracez}) is an elementary consequence of
the Cauchy-Schwartz inequality.
\QED

\medskip

We now estimate the second term of $Q_0$:
since
$$  2 |\nu'_n\, \alpha_x(x,z_n)|
  \le (1-\epsilon) (\nu'_n)^2 + {1\over 1-\epsilon} \alpha_x^2(x,z_n), $$
for any $\epsilon>0$, we have
$$  \int_{-L}^L \left(  \nu'_n - \alpha_x(x,z_n) \right)^2\, dx
   \ge  \int_{-L}^L \epsilon (\nu'_n)^2 \, dx
      -  \int_{-L}^L {\epsilon\over 1-\epsilon} \alpha_x^2(x,z_n)\, dx.
$$
We now use (\ref{tracex}) and the elliptic estimate (\ref{H1})--(\ref{C_0})
to estimate
\beann
&&   \iint (\curl\alpha)^2\, dx\, dz
              -p\sum_{n=0}^N \int_{-L}^L {\epsilon\over 1-\epsilon} 
                  \alpha_x^2(x,z_n)\, dx     \\
 &&\qquad  \ge    \left(\frac12 \left( 1 + {4\over \pi^2}
     {L^2 N^2 p^2\over 4L^2 + N^2 p^2}\right)^{-2}
       - {2\epsilon\over (1-\epsilon)}\right)
              \| \vec\alpha\|_{H^1}^2 \\
   && \qquad \ge 
   \left\{ \frac12 \left( 1 + {4\over \pi^2}L^2\right)^{-2}
     - {2\epsilon\over (1-\epsilon)}\right)
              \| \vec\alpha\|_{H^1}^2,
\eeann
uniformly in $N$.
Now choose $\epsilon=
  \frac14\min\{1,\frac12 \left( 1 + {4\over \pi^2}L^2\right)^{-2}\}$.
From the Poincar\'e inequality (for $H^1$ functions on $[-L,L]$
with vanishing mean) we may then obtain the lower bound
\beann
&&   p\sum_{n=0}^N \int_{-L}^L (\nu'_n - \alpha_x(x,z_n))^2 \, dx
          + \iint_D (\curl\vec\alpha)^2\, dx\, dz
             \\
    && \qquad \ge
      \epsilon p\sum_{n=0}^N \int_{-L}^L  (\nu'_n)^2\, dx
        + \frac14 \min\left\{ 1,\left( 1 + {4\over \pi^2}L^2\right)^{-2}\right\}
                 \|\vec\alpha\|_{H^1}^2  \\
    && \qquad \ge
      \epsilon\left( 1 + {4\over \pi^2}L^2\right)^{-1}
         p\sum_{n=0}^N \|\nu_n\|_{H^1}^2
            + \frac14 \min\left\{ 1,\left( 1 + {4\over \pi^2}L^2\right)^{-2}\right\}
                 \|\vec\alpha\|_{H^1}^2 \\
    && \qquad \ge \frac14\min\left\{ 1, 
        \left( 1 + {4\over \pi^2}L^2\right)^{-3} \right\}
          \, \left[   p\sum_{n=0}^N \|\nu_n\|_{H^1}^2
            +  \|\vec\alpha\|_{H^1}^2  \right].
\eeann
We then obtain
\be\label{lambda}  \lambda\ge {1\over 4\kappa^2} 
         \min\left\{ 1, 
        \left( 1 + {4\over \pi^2}L^2\right)^{-3} \right\}. 
\ee
Note in particular that $\lambda$ is bounded away from zero
uniformly in $N$, but the bound deteriorates as either
$L\to\infty$ or $\kappa\to\infty$.  (See Remark~\ref{Ltoinfty}.)

\medskip

\noindent {\bf Step 4:}\ An upper bound on the linearization.
We define $Q_r$ to be the quadratic form representing the
second variation of energy around the solution $s+w(r,s)\in\SSS_r$,
with respect to variations $\Psi=(\mu_n,\nu_n,\vec\alpha)$ in the subspace
$W$: 
\beann
Q_r(\Psi,\Psi) &:=& \frac12 \left. {d^2\over dt^2} \right|_{t=0}
  \Omr(s+w(r,s)+t\Psi)  \\
  & =&  p\sum_{n=0}^N \int_{-L}^L \left\{
          {1\over\kappa^2} (\mu'_n)^2 +
          {1\over\kappa^2} (v'_n-a_x)^2 \mu_n^2
          + (3f_n^2-1)\mu_n^2  \right.  \\
   &&\quad \left.
      +{4\over\kappa^2} (v'_n-a_x)f_n \mu_n (\nu'_n-\alpha_x(x,z_n))
        + {1\over\kappa^2} f_n^2 (\nu'_n-\alpha_x(x,z_n))^2 \right\}\, dx\\
   && \quad +{r\over 2}\, p\sum_{n=1}^N \int_{-L}^L
       \left\{   (f_n\mu_{n-1}+f_{n-1}\mu_n)\sin\Phi_{n,n-1}
                [\nu_n-\nu_{n-1}-\int_{z_{n-1}}^{z_n} \alpha_z\, dz]
                \right.\\
    &&\quad 
       + f_n f_{n-1}\cos\Phi_{n,n-1} 
           [\nu_n-\nu_{n-1}-\int_{z_{n-1}}^{z_n} \alpha_z\, dz]^2 \\
    && \qquad \left.
    + \mu_n^2 +\mu_{n-1}^2 
     -2\mu_n\mu_{n-1} \cos\Phi_{n,n-1}
              \right\}\, dx \\
     &&\quad + {1\over\kappa^2}
      \iint_B (\curl\vec\alpha)^2\, dx\, dz.
\eeann

As long as the linearization $Q_r$ remains non-singular
(as an operator on $W$), we
may invoke the Implicit Function Theorem to ensure that the
solutions to equation (\ref{Weqn}) in the form $s+w(r,s)$
determined by Lemma~\ref{Lya} remain valid.  Hence, the critical
value $r_0$ which determines the radius of convergence of the
expansion is bounded below by the smallest value of $r$ for which
$Q_r$ admits zero as an eigenvalue.  If $Q_r$ degenerates
at $r_*$, then there exists a test vector $\Psi\in W$,
$\|\Psi\|=1$, such that $Q_{r_*}(\Psi,\Psi)=0$.  Using the
lower bound on $Q_0$, we have:
\be\label{setup}  \lambda\le Q_0(\Psi,\Psi) = Q_{r_*}(\Psi,\Psi) +
   \left( Q_0-Q_{r_*}\right)(\Psi,\Psi)
      = \left( Q_0-Q_{r_*}\right)(\Psi,\Psi).  
\ee
An estimate of the difference between the linearization
at $r=r_*$ and at $r=0$ in terms of $r_*$ (and the norm of $\Psi$)
will then yield a lower bound on the critical value $r_0\ge r_*$.
This estimate follows from the sup-norm bounds on each 
$u_n$, $(v'_n-a_x(\cdot,z_n))$
obtained in (\ref{bestu}) and (\ref{bestj}):
\beann
|Q_{r_*}(\Psi,\Psi)- Q_0(\Psi,\Psi)|
 &\le &   p\sum_{n=0}^N \int_{-L}^L \left\{
          {1\over\kappa^2} (v'_n-a_x)^2 \mu_n^2
          + 3|f_n^2-1|\mu_n^2  \right.  \\
   &&\quad 
      +{4\over\kappa^2} |v'_n-a_x| |\mu_n| |\nu'_n-\alpha_x(x,z_n)| \\
   &&\quad \left.
      + {1\over\kappa^2} |f_n^2-1| 
            (\nu'_n-\alpha_x(x,z_n))^2 \right\}\, dx\\
   && \quad +{r_*\over 2}\, p\sum_{n=1}^N \int_{-L}^L
       \left\{   (|\mu_{n-1}|+|\mu_n|)
                |\nu_n-\nu_{n-1}-\int_{z_{n-1}}^{z_n} \alpha_z\, dz|
                \right.\\
    &&\quad \left. + \mu_n^2 +\mu_{n-1}^2 
     -2|\mu_n\mu_{n-1}|
       + [\nu_n-\nu_{n-1}-\int_{z_{n-1}}^{z_n} \alpha_z\, dz]^2
              \right\}\, dx \\
   &\le &  cr_* \,[1+ K(1+r_*\kappa^2 K)]\,  \|\Psi\|^2,
\eeann
where 
$$  K:= {1\over Hp} (1+r_*L^2\kappa^2)(1+r_*L).  $$
Inserting the above estimate 
and (\ref{lambda}) into (\ref{setup}) produces a
lower bound on the first point of degeneracy $r_*$,
\be\label{radius}
r_* \, [1+K(1+r_*\kappa^2 K)] \ge \lambda \ge 
    {1\over 4\kappa^2} 
         \min\left\{ 1, 
        \left( 1 + {4\over \pi^2}L^2\right)^{-3} \right\}.
\ee
We observe that the left-hand term is monotone increasing
in $r_*$, and therefore implies a lower bound of the
form
$$   r_0 \ge r_* \ge R(L,\kappa,H),  $$
where the constant $R$ can be chosen uniformly in $N$, but
decreases with increasing $L$ and $\kappa$ and increases with
increasing $H$.

\bigskip

\noindent {\bf Step 5:}\ Conclusion.
The above estimates have all been based on the (unsatisfactory)
initial hypothesis $|r|<r_1(N,L, \kappa, H)$ chosen small enough to
ensure (\ref{fnotzero}).  However, we observe that the only
role of this hypothesis is to obtain (\ref{fnotzero}),
and the estimates obtained thereafter are valid whenever
(\ref{fnotzero}) holds.  Since
$\|u_n\|_\infty$ is  a continuous function of $r$,
the estimates obtained for $u_n,v_n,\vec a$ and $r_*$ above
persist as we increase $r$, either until (\ref{fnotzero})
is first violated or until we reach $r=r_*$.  
From (\ref{bestu}), if condition (\ref{fnotzero}) is violated then
$r$ satisfies
$$  C_u r[1+r\kappa^2 K^2] \ge \frac12,  $$
a condition on $r$ which also defines a lower bound
which is independent of $N$, decreasing with
increasing $L,\kappa$, and increasing  as $H$ increases (just
 as for  $R(L,\kappa,H)$
above.) In either case, the perturbation argument remains valid
for an interval  $r\in [0,r_0)$ with $r_0$ bounded below
by a constant which is
independent of $N$, with the
dependences on $L,\kappa,H$ claimed in the
statement of (iv).  This concludes
the proof of Theorem~\ref{Lya}.
\QED

\bigskip

\begin{rem}\label{Ltoinfty}\rm
We note that the lower bound on $r_*$ obtained above improves
with smaller $L$, $\kappa$ and larger $H$.  Of course this is
only a one-sided bound, and we cannot be sure that the solutions
obtained by the degenerate perturbation method really do
cease to exist if $L$ or $\kappa$ are too large.  Indeed, 
solutions with a similar form may persist beyond the range
of validity of Lemma~\ref{Lya}, and the estimates which gave
the lower bound (\ref{radius}) are certainly not sharp.
However, we can show that the smallest eigenvalue $\lambda(N,L,\kappa)$
of the linearization at $r=0$  does tend to zero as either
$L,\kappa\to\infty$.  From the Implicit Function Theorem
we may then infer that $w(r,s)$ grows rapidly near $r=0$,
which in turn suggests that the interval of validity is
diminished as $L$, $\kappa$ increase.
\end{rem}
\begin{lem}\label{lambdatozero}
For any fixed $N$,
$$   \lambda(N,L,\kappa) \le{9\over  2\kappa^2    p^2 L^2} \to 0  $$
as $L,\kappa\to\infty$.
\end{lem}
\Pf
As before, set $M=Np$.  Let $\xi\in H^2\cap H^1_0(D)$ be the
solution to the Dirichlet problem
$$  \Delta\xi =1, \quad
   \xi|_{\partial D}=0,  $$
and define $\vec \alpha=(\partial_z\xi, -\partial_x\xi)$.
Note that $\xi$ is even in $x$, and
$$   {1\over\kappa^2}\iint_D (\curl\vec \alpha)^2\, dx\, dz
      ={2LM\over\kappa^2}. 
$$

We next choose $\nu_n$ with
$$  {d\over dx}
            \nu_n(x)=\alpha_x(x,z_n), \qquad \int_{-L}^L \nu_n(x)\, dx =0.  $$
Since $\alpha_x(x,z)$ is even in $x$ we have
$$  \nu_n(x)=\int_0^x \alpha_x(x',z_n)\, dx',  $$
and $\nu_n$ is odd.  Applying Stokes' Theorem to the
rectangle $[-L,L]\times [z_{n-1},z_n]$,
\beann
2px 
    &=& \int_{z_{n-1}}^{z_n}\int_{-x}^x
      \curl\vec \alpha \, dx \, dz \\
&=& \int_{-x}^x \left(\alpha_x(x,z_n) - \alpha_x(x,z_{n-1})\right)\, dx
   + \int_{z_{n-1}}^{z_n} \left( \alpha_z(-x,z) - \alpha_z(x,z)\right)\, dz \\
&=& 2\left(  \nu_n(x)-\nu_{n-1}(x) -\int_{z_{n-1}}^{z_n} \alpha_z(x,z)\,
     dz\right).
\eeann
In particular,
\beann
2\sum_{n=0}^N \|\nu_n\|_2^2  &\ge &
\sum_{n=1}^N \left[\|\nu_n\|_2^2 + \|\nu_{n-1}\|_2^2\right]  \\
  &\ge & {2 N p^2\over 9} L^3 - 
         \sum_{n=1}^N\int_{-L}^L
     \left( \int_{z_{n-1}}^{z_n} \alpha_z(x,z)\, dz\right)^2\, dx\\ 
  &\ge &  {2 N p^2\over 9} L^3 -
      p \iint_D [\alpha_z(x,z)]^2\, dz\, dx.
\eeann

Taking $\mu_n=0$, $\Psi=(\mu_n,\nu_n,\vec \alpha)$,
we have 
$$
\|\Psi\|_E^2 \ge  p\sum_{n=0}^N \|\nu_n\|^2_2 + \|\vec\alpha\|_2^2
   \ge { M p^2\over 9} L^3,
$$
assuming (as usual) that $0<p\le 1$.  We may then conclude that:
\beann
\lambda(N,L,\kappa) &\le &
    { Q_0(\Psi,\Psi)\over \|\Psi\|_E^2} \\
&\le &
   { {1\over\kappa^2}\iint_D (\curl\vec\alpha)^2\, dx\, dz
      \over  {M p^2\over 9} L^3 }  \\
& \le &  {9\over  2\kappa^2    p^2 L^2}.
\eeann
\QED

\begin{figure}\label{vplane}
\centering
\includegraphics{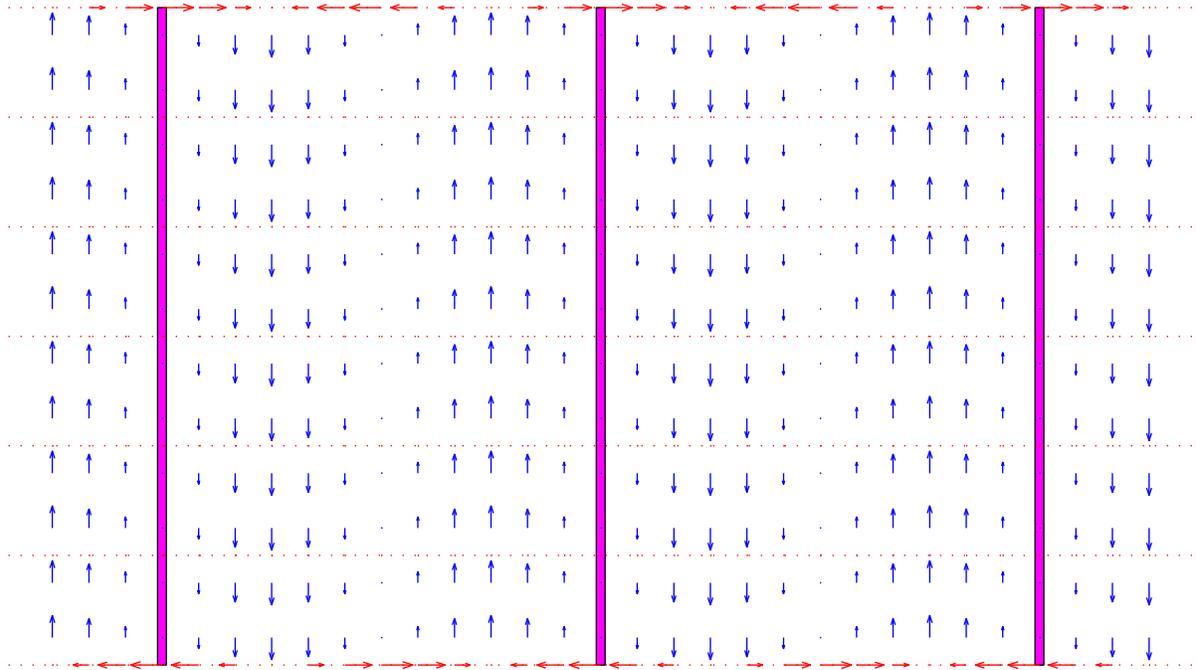}
\caption{{ \small\sl Vortex planes, for a sample with a 
finite number of superconducting planes.  The dotted
horizontal lines represent the superconducting planes,
and horizontal arrows indicate the in-plane currents $j_x^{(n)}$,
which vanish to order $r$ in the interior of the sample.
The vertical arrows indicate the Josephson currents $j_z^{(n)}$
in each gap.  Except for an edge
effect near the top and bottom of the sample, the magnetic field $h(x,z)$
and Josephson currents $j_z$ are $z$-independent.  The
vortices correspond to local maxima of $h$, and lie along planes
$x=$constant, indicated here by the dark bands.}}
\end{figure}


\begin{thebibliography}{11}
%
\bibitem[AlBeBr 00]{ABB2} S. Alama, A.J. Berlinsky, and L. Bronsard,
{\sl Periodic vortex lattices for the Lawrence--Doniach model
of layered superconductors in a parallel field,} preprint 2000,
available on the preprint archive http://xxx.lanl.gov,
math.AP/0010111.
%
\bibitem[AmBa 98]{AB} A. Ambrosetti and M. Badiale,
{\sl Homoclinics: Poincar\'e-Melnikov type results via a variational
approach,}
 Ann. Inst. H. Poincar\'e Anal. Non Lin\'eaire, vol. 15
(1998), pp. 233--252.
%
\bibitem[AmCzE 87]{ACE} A. Ambrosetti, V. Coti-Zelati, and I. Ekeland,
{\sl Symmetry breaking in Hamiltonian systems,}
J. Diff. Eq., vol. 67 (1987), pp 165--184.
%
\bibitem[BLR 95]{BLR} A. Bahri, Y. Li, O. Rey, {\sl
On a variational problem with lack of compactness: the topological effect of
the critical points at infinity,}
Calc. Var.
Partial Differential Equations, vol. 3 (1995), pp. 67--93.
%
\bibitem[BBH 94]{BBH} F. Bethuel, H. Brezis, and F. H\'elein,
{\sl Ginzburg--Landau Vortices.}  Birkhauser, Boston, 1994.
%
\bibitem[BR 95]{BR} F. Bethuel, T. Rivi\'ere, {\sl
Vortices for a variational problem related to superconductivity,}
Ann. Inst. H.
Poincar\'e Anal. Non Lin\'eaire, vol. 12  (1995), pp. 243--303.
%
\bibitem[Bu 73]{Bul} L. Bulaevskii, {\sl Magnetic properties of layered
superconductors with weak interaction between the layers,}
Sov. Phys. JETP, vol. 37 (1973), pp 1133--1136.
%
\bibitem[BuCm 91]{BC} L. Bulaevskii and J. Clem,
{\sl Vortex lattice of highly anisotropic layered superconductors
in strong, parallel magnetic fields,} Phys. Rev. B44 (1991), 
pp 10234--10238.
%
\bibitem[CDG 95]{CDG} S.  Chapman,  Q. Du, and M. Gunzburger,
{\sl On the Lawrence-Doniach and anisotropic Ginzburg-Landau models for
layered superconductors,} SIAM J. Appl. Math., vol. 55 (1995), pp.
156--174.
%
\bibitem[CmCo 90]{CC} J. Clem and M. Coffey,
{\sl Viscous flux motion in a Josephson--coupled layer model
of high--$T_c$ superconductors,}
Phys. Rev. vol. B42 (1990), pp. 6209--6216.
%
\bibitem[DF 97]{DF} M. Del Pino and P. Felmer, {\sl
Local minimizers for the Ginzburg-Landau energy,}
Math. Z., vol. 225 (1997), pp. 671--684.
%
\bibitem[GiP 99]{GP}  T. Giorgi and D. Phillips, {\sl
The breakdown of superconductivity due to strong
fields for the Ginzburg--Landau model,}
SIAM J. Math. Anal., vol. 30 (1999), pp. 341--359.
%
\bibitem[Gr 85]{Grisvard} Pierre Grisvard, ``Elliptic Problems in
Nonsmooth Domains,'' Pitman Advanced Publishing Program,
Boston, 1985.
%
\bibitem[Gu 96]{Gui} C. Gui, {\sl
Multipeak solutions for a semilinear Neumann problem,}
Duke Math. J., vol. 84 (1996), pp. 739--769. 
%
\bibitem[I 92]{Iye} Y. Iye, {\sl How Anisotropic Are the 
Cuprate High $T_c$ Superconductors?} Comments Cond. Mat. Phys.,
vol. 16 (1992), pp. 89--111.
%
\bibitem[KAVB 90]{Kes} P. Kes, J. Aarts, V. Vinokur, and
C. van der Beek, {\sl Dissipation in Highly Anisotropic Superconductors,}
Phys. Rev. Lett. vol. 64 (1990), pp 1063--1066.
%
\bibitem[K 99]{K} S. Kuplevakhsky, {\sl Microscopic theory of
weakly couple superconducting multilayers in an external
magnetic field,} preprint {\tt cond-mat/9812277}.
%
\bibitem[LaDo 71]{LD} W. Lawrence and S. Doniach,
{\it Proceedings of the Twelfth
International Conference on Low Temperature Physics,}
E. Kanda (ed.), Academic Press of Japan, Kyoto, 1971, p. 361.
%
\bibitem[LN 98]{LN} Y. Li and L. Nirenberg, {\sl 
The Dirichlet problem for singularly perturbed elliptic equations,}
Comm. Pure Appl. Math., vol. 51 (1998), pp. 1445--1490.
%
\bibitem[LL 97]{LiebLoss} E. Lieb and M. Loss
{\sl Analysis}. Graduate Studies in Mathematics, vol. 14. 
American Mathematical
Society, Providence, RI, 1997.
%
\bibitem[RuS 98]{RS}  J. Rubinstein and M. Schatzman,
{\sl Asymptotics for thin superconducting rings,}  J. Math. Pures Appl., 
s\'erie 9,
vol. 77  (1998), pp. 801--820.
%
\bibitem[Re 91]{Rey} O. Rey, {\sl 
Blow-up points of solutions to elliptic equations with limiting
nonlinearity,}  Differential Integral Equations, vol. 4 (1991), pp.
1155--1167.  
%
\bibitem[Th 90]{Th} S. Theorodakis, {\sl Theory of
vortices in weakly-Josephson-coupled layered superconductors,}
Phys. Rev. B42 (1990), pp 10172--10177.
%
\bibitem[T 96]{Tinkham} M. Tinkham, ``Introduction to
Superconductivity,'' 2nd edition, Mc Graw-Hill, New York, 1996.
%
\bibitem[W 98]{W} J. Wei, {\sl
On the interior spike solutions for some singular perturbation problems,}
Proc. Roy. Soc. Edinburgh, Sect. A, vol. 128 (1998), pp. 849--874. 


\end{thebibliography}
\end{document}